\documentclass[a4paper,11pt]{amsart}
\usepackage{fancyhdr}
\usepackage{amsmath}
\usepackage{dsfont}
\usepackage{hyperref}
\usepackage[mathscr]{eucal}
\usepackage[cp1251]{inputenc}
\usepackage[english]{babel}
\usepackage{cite,enumerate,float,indentfirst}
\usepackage{graphicx}
\usepackage{xcolor}
\usepackage{latexsym,a4,mathrsfs,amsthm,amsmath,amssymb,url}
\usepackage{amsfonts}
\usepackage{amssymb}
\usepackage{epstopdf}

\numberwithin{equation}{section}
\newtheorem{theorem}{Theorem}[section]
\newtheorem{lemma}[theorem]{Lemma}
\newtheorem{statement}[theorem]{Statement}
\newtheorem{definition}[theorem]{Definition}
\newtheorem{corollary}[theorem]{Corollary}
\theoremstyle{remark}
\newtheorem*{remark}{Remark}

\newcommand{\C}{\mathbb{C}}

\newcommand{\Cond}{{\bf (C0)}}
\newcommand{\UI}{{\bf (UI)}}
\newcommand{\Lind}{{\bf (L)}}

\newcommand{\X}{{\bf X}}

\newcommand{\M}{{\bf M}}

\newcommand{\HH}{{\bf H}}

\newcommand{\I}{{\bf I}}
\newcommand{\RR}{{\bf R}}
\newcommand{\V}{{\bf V}}
\newcommand{\J}{{\bf J}}
\newcommand{\W}{{\bf W}}

\newcommand{\ee}{{\bf e}}

\DeclareMathOperator{\Tr}{Tr}

\DeclareMathOperator{\E}{\mathbb{E}}
\DeclareMathOperator{\Pb}{\mathbb{P}}

\DeclareMathOperator{\imag}{Im}

\begin{document}

\vspace{1in}

\title[Product of random matrices]{\bf On One Generalization of the Elliptic Law for Random Matrices}

%\vspace{2in}

\author[F. G{\"o}tze]{F. G{\"o}tze}
\address{F. G{\"o}tze\\
 Faculty of Mathematics\\
 Bielefeld University \\
 Bielefeld, Germany
}
\email{goetze@math.uni-bielefeld.de}

\author[A. Naumov]{A. Naumov}
\address{A. Naumov\\
 Faculty of Computational Mathematics and Cybernetics\\
 Moscow State University \\
 Moscow, Russia
 }
\email{naumovne@gmail.com, anaumov@math.uni-bielefeld.de}

\author[A. Tikhomirov]{A. Tikhomirov}
\address{A. Tikhomirov\\
 Department of Mathematics\\
 Komi Research Center of Ural Division of RAS \\
 Syktyvkar, Russia
 }
\email{tichomir@math.uni-bielefeld.de}
\thanks{All authors are supported by CRC 701 ``Spectral Structures and Topological
Methods in Mathematics'', Bielefeld. A.~Tikhomirov and A.~Naumov are partially supported by RFBR, grant  N 14-01-00500 ``Limit theorems for random matrices and and their applications''. A.~Tikhomirov are supported by Program of Fundamental Research Ural Division of RAS, Project N 12-P-1-1013}


\date{\today}

\begin{abstract}
We consider the products of $m\ge 2$ independent large real random matrices with independent vectors
$(X_{jk}^{(q)},X_{kj}^{(q)})$ of entries. The entries $X_{jk}^{(q)},X_{kj}^{(q)}$ are correlated with $\rho=\E X_{jk}^{(q)}X_{kj}^{(q)}$. The limit distribution of the empirical spectral distribution of the eigenvalues of such products doesn't depend on $\rho$ and equals to the distribution of $m$th power of the random variable uniformly distributed on the unit disc.
\end{abstract}

\maketitle

\section{Introduction}
Let $m\ge 1$ be a fixed integer and $\X^{(q)} = n^{-1/2} \{X_{jk}^{(q)}\}_{j,k=1}^n, q = 1, ... , m$, be independent random matrices with real entries. We suppose that the random variables $X_{j,k}^{(q)}$, $1\le j,k\le n, q = 1, ... ,m$, are defined on a common probability space $(\Omega,\mathcal F,\Pb)$ and satisfy the following conditions \Cond:\\
a) random vectors $(X_{jk}^{(q)}, X_{kj}^{(q)})$ are mutually independent for $1 \le j < k \le n$;\\
b) for any $1 \le j \le k \le n$
$$
\E X_{j k}^{(q)} = 0 \text{ and } \E (X_{j k}^{(q)})^2 = 1;
$$
c) for any $1 \le j < k \le n$
$$
\E ( X_{j k}^{(q)} X_{k j}^{(q)} ) = \rho, |\rho| \le 1;
$$
d) diagonal entries and off-diagonal entries are independent.

We say that the random variables $X_{j,k}^{(q)}$, $1\le j,k\le n, q = 1,...,m$, satisfy the condition $\UI$ if the squares of $X_{jk}^{(q)}$'s are uniformly integrable , i.e.
\begin{equation}\label{unif0}
\max_{q,j,k}\E|X_{jk}^{(q)}|^2\mathbb I{\{|X_{jk}^{(q)}|>M\}} \rightarrow 0 \quad \text{as}\quad M \rightarrow \infty.
\end{equation}
Here and in what follows $\mathbb I\{B\}$ denotes the indicator of the event $B$.

The random variables $X^{(q)}_{jk}$ may depend on $n$, but for simplicity we shall not make this explicit in our  notations. Denote by $\lambda_1, ..., \lambda_n$ the eigenvalues of the matrix $\mathbf W:= \prod_{q=1}^m\mathbf X^{(q)}$ and define the empirical spectral measure of the eigenvalues by
$$
\mu_n(B) = \frac{1}{n} \# \{1 \le i \le n: \lambda_i \in B \}, \quad B \in \mathcal{B}(\C),
$$
where $\mathcal{B}(\C)$ is a Borel $\sigma$-algebra of $\C$.

We say that the sequence of random probability measures $m_n(\cdot)$
converges weakly in probability to the probability measure $m(\cdot)$ if for all continues and bounded functions $f: \C \rightarrow \C$ and all $\varepsilon > 0$
$$
\lim_{n \rightarrow \infty}\Pb \left ( \left | \int_\C f(x) m_{n}(dz) - \int_\C f(x)m(dz) \right | > \varepsilon \right ) = 0.
$$
We denote a weak convergence by the symbol $\xrightarrow{weak}$.

A fundamental problem in the theory of random matrices is to determine the limiting distribution of $\mu_n$
as the size of the random matrix tends to infinity. The following theorem gives the solution of this problem for the matrices which satisfy the conditions $\Cond$ and $\UI$.

\begin{theorem} \label{th:main}
Let $m\ge 2$ and $\X^{(q)} = n^{-1/2}\{X_{jk}^{(q)}\}_{j,k=1}^n, q = 1, ... , m,$ be independent random matrices such that the random variables $X_{jk}^{(q)}, j,k=1,...,n, q = 1, ..., m,$ satisfy the conditions $\Cond$ and $\UI$. Assume that $|\rho| < 1$. Then $\mu_n \xrightarrow{weak} \mu$ in probability, and $\mu$ has the density $g$:
$$
g(x, y) = \frac1{\pi m(x^2+y^2)^{\frac{m-1}m}}\mathbb I\{x^2+y^2\le1\},
$$
which doesn't depend on $\rho$.
\end{theorem}
\begin{remark}
Theorem~\ref{th:main} was announced in the talk of F. G{\"o}tze ``Spectral Distribution of Random Matrices and Free Probability'', Advanced School and Workshop on Random Matrices and Growth Models, Trieste, Italy. Recently O'Rourke, Renfrew, Soshnikov and Vu, see~\cite{RRSV2014} , proved the result of Theorem~\ref{th:main} under additional assumptions on the moments of $X_{jk}^{(q)}$.
\end{remark}
\begin{remark}
Girko~\cite{Girko1985} showed  that  for $m=1$ under the additional assumptions that the distribution of r.v.'s $X_{jk}^{(1)}$ has a density the limit measure $\mu$ has a density of uniform distribution on the ellipse
$\mathcal E=\{(x,y): \frac{x^2}{(1-\rho)^2}+\frac{y^2}{(1+\rho)^2}\le 1\}$.  This result is called ``elliptic law''.
For Gaussian matrices the elliptic law was proved in~\cite{SomStein1988}. The elliptic law without assumption on the density of distribution of  entries $X_{jk}$ was proved by Naumov in~\cite{Nau2012}. Nguyen and O'Rourke in~\cite{ORourkeNguyen2012} and G\"otze, Naumov, Tikhomirov in~\cite{GotNauTikh2013sing} extended the elliptic law on the case when $X_{jk}^{(1)}$'s have only finite second moment and non-identical distribution.
\end{remark}
\begin{remark}
For $m=1$ and $\rho=0$ we have the circular law, i.e. the limit distribution $\mu$ is uniform distribution on the unit disc.
The circular law was first proved by Ginibre in\cite{Ginibre1965} for matrices with independent standard complex Gaussian entries.  Girko in  \cite{Girko1984} have considered the general case under assumption that the distributions of entries have bounded densities and  the fourth moments of entries are finite.  Z. Bai (see \cite{BaiSilv2010}) rely on the fruitful Girko's ideas gave a correct proof of the circular law under the same assumptions. G\"otze and Tikhomirov in \cite{GotTikh2007} have proved the circular law without assumption on the density of entries, but assuming the sub-Gaussian distributions of r.v.'s $X_{jk}^{(1)}$. Later Pan and Zhou in \cite{PanZhou2007}
proved the circular law assuming that $\E|X_{jk}^{(1)}|^4 < \infty$. G\"otze and Tikhomirov in \cite{GotTikh2010} proved the circular law assuming the logarithmic second moments ($\E|X_{jk}^{(1)}|^2|\log|X_{jk}^{(1)}||^\alpha<\infty$ with some $\alpha$ sufficiently large). And finally Tao and Vu in \cite{TaoVu2010} proved the Circular law for i.i.d. case under the assumption on the second moments only.
\end{remark}
\begin{remark}
In the case $\rho=0$ and $X_{jk}^{(q)}$ and $X_{kj}^{(q)}$ are independent for $1\le j<k\le n$, Theorem \ref{th:main} was proved by G\"otze and Tikhomirov in \cite{GotTikh2010}. See also the result of O'Rourke and Soshnikov~\cite{RS2011}.
\end{remark}

%%%%%%%%%%%%%%%%%%%%%%%%%%%%%%%%%%%%%%%%%%%%%%%%%%%%%%%%%%%%%%%%%%%%%%%%%%%%%%%%%%%%%%%%%%
\subsection{Proof of the elliptic law}
In the  following  we shall give  the proof of Theorem~\ref{th:main}. We shall use the logarithmic potential
approach first suggested for the proof of the circular law by G\"otze and Tikhomirov in \cite{GotTikh2007}.
This approach was developed in many papers (see, for instance \cite{GotTikh2010}, \cite{GotTikh2011} and \cite{BordCh2011}).
We define the logarithmic potential of the empirical spectral measure of the matrix $\W$ by the formula
\begin{equation*}
 U_n(z)=-\int_{\C}\ln|w-z| \mu_n(dw)
\end{equation*}
and will prove that
\begin{equation*}
 \lim_{n\to\infty} U_n(z)=U(z):=-\int_{\C}\ln|w-z|\mu(dw).
\end{equation*}
Let us denote by $s_1 \geq s_2 \geq ... \geq s_n$ the singular values of $\W - z \I$ and introduce the empirical spectral measure $\nu_n(\cdot, z)$ of squares of singular values. We can rewrite the logarithmic potential of $\mu_n$ via the logarithmic moments of measure $\nu_{n}$ by
\begin{align*} \label{eq:log. potential}
&U_{\mu_n} (z) = -\int_\C \ln|z - w| \mu_n(dw) = -\frac{1}{n} \ln \left |\det \left (\W - z\I \right )\right| \\
&= - \frac{1}{2 n} \ln \det \left (\W-z\I \right )^{*}\left (\W-z\I \right ) =
-\frac{1}{2}\int_0^\infty \ln x \nu_{n}(dx).
\end{align*}
This allows us to consider the Hermitian matrices $(\W-z\I)^{*}(\W-z\I)$ instead of $\W$. To prove Theorem~\ref{th:main} we need the following lemma.
\begin{lemma} \label{l:Girko}
Suppose that for a.a. $z \in \mathbb C$ there exists a probability measure $\nu_z$ on $[0, \infty)$ such that\\
a) $\nu_n \xrightarrow{weak} \nu_z$ as $n \rightarrow \infty$ in probability\\
b) $\ln$ is uniformly integrable in probability with respect to $\{\nu_n\}_{n \geq 1}$.

Then there exists a probability measure $\mu$ such that\\
a) $\mu_{n} \xrightarrow{weak} \mu$ as $n \rightarrow \infty$ in probability \\
b) for a.a. $z \in \C$
$$
U_{\mu}(z) = - \int_{0}^\infty \ln x \nu_z(dx).
$$
\end{lemma}
\begin{proof}
See \cite{BordCh2011}[Lemma~4.3] for the proof.
\end{proof}
\begin{proof}[Proof of Theorem~\ref{th:main}]
From Lemma~\ref{l:Girko} it follows that to prove Theorem~\ref{th:main} it is enough to check conditions a) and b) and show that $\nu_z$ determines the logarithmic potential of the measure $\mu$. In Theorem~\ref{shift} we find the limit distribution of singular values of the shifted matrix $\mathbf W(z)=\mathbf W-z\mathbf I$ (Section \ref{sec2}). The solution of this problem is divided into several steps. We make symmetrization of one-sided distribution functions. Then we reduce the problem to the case of truncated random variables. Next we show that the limit of empirical distribution of singular values of product of  matrices with truncated random variables is the same as one of the product of matrices with Gaussian entries. Finally, we  show that the limit of expected distributions of singular values of matrices with Gaussian entries exists and its Stieltjes transform $s(z)$ satisfies the following system of  equations
\begin{align}
&1+w s(\alpha,z)+(-1)^{m+1}w^{m}s(\alpha,z)^{m+1}=0,\notag\\
&(w-\alpha)^2+(w-\alpha)-4|z|^2s(\alpha,z)=0.\notag
\end{align}
From the paper~\cite{GotTikh2011} we know that the measure with the Stieltjes transform $s(z)$ which satisfies this system of  equations determines the logarithmic potential of the measure $\mu$.

In Section \ref{singular}, Lemma~\ref{l: log uniform integr} we show that $\ln(\cdot)$ is uniformly integrable in probability with respect to $\{\nu_n\}_{n \geq 1}$.
\end{proof}

By $C$ (with an index or without it) we shall denote generic absolute constants, whereas $C(\,\cdot\,,\,\cdot\,)$ will denote positive constants depending on arguments. {For any matrix $\mathbf A$ we shall denote by $\|\mathbf A\|_2$ the Frobenius norm of matrix $\mathbf A$ ($\|\mathbf A\|_2^2=\Tr\mathbf A\mathbf A^*$) and by $\|\mathbf A\|$ we shall denote the operator norm of matrix $\mathbf A$ ($\|\mathbf A\|=\sup_{\mathbf x: \|\mathbf x\|=1}\|\mathbf A\mathbf x\|$). Here and in the what follows $\mathbf A^*$ denotes the adjoined (transposed and complex conjugate)  matrix $\mathbf A$}

%%%%%%%%%%%%%%%%%%%%%%%%%%%%%%%%%%%%%%%%%%%%%%%%%%%%%%%%%%%%%%%%%%%%%%%%%%%%%%%%%%%%%%%%%%

\section{The limit distribution for singular values distribution of shifted matrices}\label{sec2}
In this Section we prove that there exists the limit distribution for the empirical spectral distribution of the matrices $\mathbf W-z\mathbf I$. Let $s_1\ge\ldots\ge s_n$ denote the singular values of the matrix $\mathbf W-z\mathbf I$. By $\mathcal G_n(x,z)$ we denote the empirical spectral distribution function of the matrix $(\mathbf W-z\mathbf I)(\mathbf W-z\mathbf I)^*$ (the distribution function of the uniform distribution on the squared singular values of the matrix $\mathbf W-z\mathbf I$). This distribution function corresponds to the measure $\nu_n(\cdot, z)$ introduced in the previous section. Let $G_n(x,z):=\E\mathcal G_n(x,z)$.

We say the entries $X_{j,k}^{(q)}$, $1\le j,k\le n, q = 1, ..., m$, of the matrices $\X^{(q)}$ satisfy Lindeberg's condition $\Lind$ if
$$
\text{for all $\tau>0$ } \quad L_n(\tau):=\max_{q = 1, ... , n}\frac{1}{n^2}\sum_{i,j=1}^n \E X_{ij}^2 \mathbb I(|X_{ij}| \geq \tau \sqrt n) \rightarrow 0 \text{ as } n \rightarrow \infty.
$$
It easy to see that $\UI \Rightarrow \Lind$

We prove the following Theorem
\begin{theorem}\label{shift}
Let $X_{jk}^{(q)}$'s satisfy the conditions $\Cond$ and $\UI$. Then there exists a distribution function $G(x,z)$ such that:\\
1) $G_n(x,z)\to G(x,z)$ as $n\to \infty$;\\
2) the Stieltjes transform  $s(\alpha,z)$ of the distribution function $G(x,z)$, defined by the equality
$s(\alpha,z):=\int\frac1{x-\alpha}dG(x,z)$, satisfies the following system of equations:
\begin{align*}
&1+w s(\alpha,z)+(-1)^{m+1}w^{m}s(\alpha,z)^{m+1}=0\\
&(w-\alpha)^2+(w-\alpha)-4|z|^2s(\alpha,z)=0,
\end{align*}
where $\imag (w-\alpha)>0$ for $\imag\alpha>0$.
\end{theorem}

\begin{remark}
It is well-known that the distribution function with Stieltjes transform satisfying the system exists and unique. Moreover, this distribution is finitely supported and has a density. (See, for instance \cite{GotTikh2011}). In particular, if $G_n(x,z)$ convergence to $G(x,z)$ then this convergence is uniformly in $x\in\mathbb R$, i.e.
\begin{equation}\notag
 \lim_{n\to\infty}\Delta_n(z)=\sup_x| G_n(x,z)-G(x,z)|\to0.
\end{equation}
\end{remark}
\begin{remark}
By Lemma~\ref{var0} one may show that $\mathcal G_n(x,z)$ weakly converges in probability to $G(x,z)$.
\end{remark}

%%%%%%%%%%%%%%%%%%%%%%%%%%%%%%%%%%%%%%%%%%%%%%%%%%%%%%%%%%%%%%%%%%%%%%%%%%%%%%%%%%%%%%%%%%

\subsection{The proof of Theorem \ref{shift}}
As we noted before we divide the proof into several steps.
\subsubsection{Symmetrization}
We will use the following  ``symmetrization'' of one-sided distributions. Let $\xi^2$ be a positive random variable with the distribution function $F(x)$. Define $\widetilde \xi:=\varepsilon\xi$ where $\varepsilon$ denotes a Rademacher random variable with $\Pb\{\varepsilon=\pm1\}=1/2$ which is independent of $\xi$. Let $\widetilde F(x)$ denote the distribution function of $\widetilde \xi$. It satisfies the equation
\begin{equation}\label{sym}
\widetilde F(x)=1/2(1+\text{sgn} \{x\}\,F(x^2)),
\end{equation}

\begin{lemma}\label{symmetrization}
For any one-sided distribution function $F(x)$ and $G(x)$ we have
\begin{equation*}
\sup_{x\ge 0}| F(x)-G(x)|=2\sup_x|\widetilde F(x)-\widetilde G(x)|,
\end{equation*}
where $\widetilde F(x)$ ($\widetilde G(x)$) denotes the symmetrization of $F(x)$ ($G(x)$ respectively) according to (\ref{sym}).
\end{lemma}
\begin{proof}By (\ref{sym}), we have for any $x\ge 0$
\begin{align*}
F(x)=2\widetilde F(\sqrt x)-1\\
G(x)=2\widetilde G(\sqrt x)-1.
\end{align*}
This implies
\begin{equation*}
\sup_{x\ge 0}| F(x)-G(x)|=2\sup_{x\ge 0}|\widetilde F(\sqrt x)-
\widetilde G(\sqrt x)|=2\sup_{x}|\widetilde F( x)-
\widetilde G( x)|.
\end{equation*}
Thus Lemma is proved.
\end{proof}
We apply this Lemma to the distribution of the  squared
singular values of the matrix $\mathbf W-z\mathbf I$. Introduce the following matrices
\begin{align}\notag
\mathbf V&=\left(\begin{matrix}{\mathbf W\quad \mathbf O}\\{\mathbf O\quad
\mathbf W^*}\end{matrix}\right),\quad
\mathbf J(z)=\left(\begin{matrix}{\mathbf O\quad z\mathbf I}\\{ z\mathbf I\quad
\mathbf O}\end{matrix}\right),
\quad\mathbf J=\mathbf J(1),\quad\text{and}\quad
 {\mathbf V}(z)=\mathbf V\mathbf J-\mathbf J(z),\notag\\
 \mathbf R&:=(\mathbf V(z)-\alpha\mathbf I)^{-1},\notag
\end{align}
where $\mathbf I$ denotes the unit matrix of the corresponding order and $\alpha=u+iv\in \mathbb C^+$ $(v>0)$.
Note that ${\mathbf V}(z)$ is a Hermitian matrix. The eigenvalues of the matrix ${\mathbf V}(z)$ are $-s_1,\ldots,-s_n,s_n,\ldots,s_1$. Note that the symmetrization of the distribution function $\mathcal
G_n(x,z)$ is a function $\widetilde{\mathcal G}_n(x,z)$ which is the empirical distribution function  of the   eigenvalues of  the matrix ${\mathbf V}(z)$. According to Lemma \ref{symmetrization}, we get
\begin{equation*}
\Delta_n(z):=\sup_x| \mathcal G_n(x,z)-G(x,z)|=
2\sup_x|\widetilde {\mathcal G}_n(x,z)-\widetilde G(x,z)|=:2\widetilde \Delta_n(z).
\end{equation*}
Up to now we shall proof that $\lim_{n\to\infty}\widetilde \Delta_n(z)=0$.
In what follows we shall consider symmetrizing distribution function only.
We shall omit symbol $"\,\widetilde {\cdot}\,"$ in the corresponding notation.

%%%%%%%%%%%%%%%%%%%%%%%%%%%%%%%%%%%%%%%%%%%%%%%%%%%%%%%%%%%%%%%%%%%%%%%%%%%%%%%%%%%%%%%%%%

\subsubsection{Truncation}
We shall now modify the random matrices $\mathbf X^{(q)}$, $q=1,\ldots,m$, by truncation of its entries.  Let $\{\tau_n\}$ is a sequence such that
\begin{equation*}
\lim_{n\to\infty}L_n(\tau_n)=0
\end{equation*}
and
\begin{equation*}
\lim_{n\to\infty}\tau_n\sqrt n=\infty.
\end{equation*}
It is well-known that such sequence there exists since $\lim_{n\to\infty}L_n(\tau)=0$ for any $\tau>0$ and $L_n(\tau)$
is non-decreasing function of $\tau$.

Introduce the random variables $X^{(q,c)}_{jk}=X^{(q)}_{jk}\mathbb I (|X^{(q)}_{jk}|\le c\tau_n\sqrt n)$  and $\overline X^{(q,c)}_{jk}=X^{(q,c)}_{jk}-\E X^{(q,c)}_{jk}$. Introduce  the  matrices $\X^{(q,c)}=\frac{1}{\sqrt n} \{X^{(q,c)}_{jk}\}_{j,k=1}^n$ and ${\overline{\X}}^{(q,c)}=\frac1{\sqrt{n}}\{{\overline X}^{(q,c)}_{jk}\}_{j,k=1}^n$.  We define the corresponding matrices $\W^{(c)}, {\overline{\W}}^{(c)}, \V^{(c)}, {\overline{\V}}^{(c)}$ and $\RR^{(c)}, {\overline{\RR}}^{(c)}$ replacing $\X^{(q)}$ in the notation of $\V, \W$ and $\RR$ by $\X^{(q,c)},{\overline{\X}}^{(q,c)}$.

Denote by $s_1^{(c)}\ge\ldots\ge s_n^{(c)}$ and $\overline s_1^{(c)}\ge\ldots\ge \overline s_n^{(c)}$ -- the singular values  of the random  matrices $ {\W}^{(c)}-z\mathbf I$ and ${\overline{\W}}^{(c)}-z\mathbf I$
respectively. We define the empirical distribution functions of the matrices $ {\V}^{(c)}(z)$ and
${\overline{\V}}^{(c)}(z)$ by
\begin{align}
{\mathcal G}_n^{(c)}(x,z)&=\frac1{2n}\sum_{k=1}^n\mathbb I(s_k^{(c)} \le x)+\frac1{2n}\sum_{k=1}^n \mathbb I (-s_k^{(c)}\le x )\notag\\
\overline{\mathcal G}_n^{(c)}(x,z)&=\frac1{2n}\sum_{k=1}^n\mathbb I({\overline s_k^{(c)}}\le
x )+\frac1{2n}\sum_{k=1}^n\mathbb I({-\overline s_k^{(c)}}\le x)\notag
\end{align}
Let
$s_n(\alpha,z)$, $ s_n^{(c)}(\alpha,z)$ and $\overline s_n^{(c)}(\alpha,z)$ denote the Stieltjes transforms of
the distribution functions $G_n(x,z)$, $G_n^{(c)}(x):=\E{\mathcal G}_n^{(c)}(x,z)$ and $\overline G_n^{(c)}(x,z)=\E{{\overline{\mathcal G}}}_n^{(c)}(x,z)$ respectively.

\begin{lemma}\label{truncation}
Under the assumptions of Theorem~\ref{th:main} the following holds:  for any $\delta>0$
\begin{equation*}
\lim_{n\to\infty}|s_n(z,\alpha)-\overline s_n^{(c)}(\alpha,z)|=0
\end{equation*}
uniformly in $\alpha=u+iv$ with $v\ge \delta$.
\end{lemma}
\begin{proof}
We compare the Stieltjes transforms $s_n(\alpha,z)$, $s_n^{(c)}(\alpha,z)$ and
$\overline s_n^{(c)}(\alpha,z)$ sequentially.
First we note that
\begin{equation}\label{resolvent}
s_n(\alpha,z)=\frac1{2n}\E\Tr \mathbf R,\qquad\text{and}\qquad s_n^{(c)}(\alpha,z)=\frac1{2n}\E\Tr \mathbf R^{(c)}.
\end{equation}
Applying the resolvent equality
\begin{equation}\notag
 (\mathbf A+\mathbf B-\alpha\mathbf I)^{-1}=(\mathbf A-\alpha\mathbf I)^{-1}-(\mathbf A-\alpha\mathbf I)^{-1}
\mathbf B(\mathbf A+\mathbf B-\alpha\mathbf I)^{-1},
\end{equation}
 we get
\begin{equation}\label{resolv}
|s_n(\alpha,z)-s_n^{(c)}(\alpha,z)|\le \frac1{2n}\E|\Tr \mathbf R^{(c)}(\mathbf V-\mathbf V^{(c)})\mathbf J\mathbf
 R|.
\end{equation}
Let
\begin{equation}\notag
\mathbf H^{(\nu)}=\left(\begin{matrix}{\mathbf X^{(\nu)}\quad\quad\quad\mathbf O}
\\{\mathbf O\quad\quad{\mathbf X^{(m-\nu+1)}}^*}\end{matrix}\right)
\quad\text{and}\quad\mathbf H^{(\nu,c)}=\left(\begin{matrix}{\mathbf X^{(\nu,c)}
\quad\quad\mathbf O}\\{\mathbf O\quad{\mathbf X^{(m-\nu+1,c)}}^*}\end{matrix}\right)
\end{equation}
Introduce the  matrices
\begin{equation}\label{vab}
\mathbf V_{a,b}=\prod_{q=a}^b\mathbf H^{(q)},\quad\mathbf V_{a,b}^{(c)}=\prod_{q=a}^b\mathbf H^{(q,c)},
\end{equation}
{($\mathbf V_{a,b}=\mathbf I$ if $a>b$).}
We have
\begin{equation}\label{repr1}
\mathbf V-\mathbf V^{(c)}=\sum_{q=1}^{m}\mathbf V^{(c)}_{1,q-1}(\mathbf H^{(q)}-\mathbf H^{(q,c)})\mathbf V_{q+1,m}.
\end{equation}
{Inequalities $\max\{\|\mathbf
R\|,\,\|\mathbf R^{(c)}\|\}\le v^{-1}$, $\|\Tr \mathbf A\mathbf B|\le \|\mathbf A\|_2\,\|\mathbf B\|_2$, inequality (\ref{resolv}),
and  the  representations~\eqref{repr1}  together imply
\begin{align}\label{st1}
|s_n(\alpha,z)-s_n^{(\alpha,c)}(z)| &\le \nonumber \\
&\frac{C}{\sqrt n}\sum_{q=1}^{m-1}\E ^{\frac12}\|(\mathbf H^{(q)}-
\mathbf H^{(q,c)})\|_2^2\frac{1}{\sqrt n}\E^{\frac12}\|\mathbf V^{(c)}_{q+1,m}\mathbf R\mathbf R^{(c)}\mathbf V_{1,q-1}\|_2^2.
\end{align}
We use here that $\Tr\mathbf A\mathbf B=\Tr\mathbf B\mathbf A$ as well.
Applying well-known   inequalities for matrix norms $\|\mathbf A\mathbf B\|_2\le \|\mathbf A\|\,\|\mathbf B\|_2$
and relation $\|\mathbf A\mathbf B\|_2=\|\mathbf B\mathbf A\|_2$ together,  we get}
\begin{equation}\notag
 \E\|\mathbf V_{q+1,m}\mathbf R\mathbf R^{(c)}\mathbf V_{1,q-1}^{(c)}\|_2^2
\le\frac C{v^4}\E\|\mathbf V^{(c)}_{1,q-1}\mathbf V_{q+1,m}\|_2^2
\end{equation}
In view  of Lemma \ref{norm2}, we obtain
\begin{equation}\label{st2}
 \E\|\mathbf V_{q+1,m}\mathbf R\mathbf R^{(c)}\mathbf V_{1,q-1}^{(c)}\|_2^2\le \frac{Cn}{v^4}.
\end{equation}
Direct calculations show that, for any $q=1,\ldots,m$,
\begin{equation}\notag
 \frac1n\E\|\mathbf X^{(q)}-\mathbf X^{(q,c)}\|_2^2\le \frac C{n^2}\sum_{j=1}^{n}\sum_{k=1}^{n}\E|X^{(q)}_{jk}|^2
I_{\{|X^{(q)}_{jk}|\ge c \tau_n\sqrt n\}}\le CL_n(\tau_n).
\end{equation}
This inequality implies that
\begin{equation}\label{st3}
\max_{1\le q\le m}\E\|\mathbf H^{(q)}-\mathbf H^{(q,c)}\|_2\le C L_n(\tau_n).
\end{equation}
Inequalities \eqref{st1}, \eqref{st2} and \eqref{st3} together imply
\begin{equation*}
 |s_n(\alpha,z)- s_n^{(c)}(\alpha,z)|\le\frac {C\sqrt {L_n(\tau_n)}}{v^2}.
\end{equation*}

Furthermore, we compare the Stieltjes transforms $s_n^{(c)}(\alpha,z)$ and $\overline s_n^{(c)}(\alpha,z)$.
By definition of $X_{jk}^{(c)}$, we have
\begin{equation}\notag
|\E X_{jk}^{(q,c)}|=|\E X_{jk}^{(q)}\mathbb I\{|X_{jk}|\ge c\tau_n\sqrt n\}|\le
 \frac1{c\tau_n\sqrt n}\E |{X_{jk}^{(q)}}|^2 I{\{|X_{jk}|\ge c\tau_n\sqrt n\}}.
\end{equation}
This implies that
\begin{equation}\label{st4}
 \|\E\mathbf X^{(q,c)}\|_2^2\le \frac C{n}\sum_{j=1}^{n}\sum_{k=1}^{n}|\E X_{jk}^{(q,c)}|^2
\le
\frac {CL_n(\tau_n)}{c\tau_n^2}.
\end{equation}

Note that ${\overline{\mathbf H}}^{(q,c)}=\mathbf H^{(q,c)}-\E\mathbf H^{(q,c)}$. Similar to the inequality (\ref{st1}) we get
\begin{equation}\notag
 |s_n^{(c)}(\alpha,z)-\overline s_n^{(c)}(\alpha,z)|\le \sum_{q=1}^{m}\frac1{\sqrt n}\|\E\mathbf H^{(q,c)}\|_2\frac1{\sqrt n}
\E^{\frac12}\|\widehat{\mathbf V}^{(c)}_{q+1,m}\mathbf R^{(c)}
\widehat{\mathbf R}^{(c)}{\widehat {\mathbf V}_{1,q-1}}^{(c)}\|_2^2.
\end{equation}
Analogously to inequality (\ref{st2}), we get
\begin{equation}\label{st100}
 \E\|\widehat{\mathbf V}^{(c)}_{q+1,m}\mathbf R^{(c)}
\widehat{\mathbf R}^{(c)}\widehat {\mathbf V}_{1,q-1}^{(c)}\|_2^2\le \frac {Cn}{v^4}.
\end{equation}
By the inequality (\ref{st4}),
\begin{equation*}
 \|\E \mathbf X^{(q,c)}\|_2\le \frac{C\sqrt {L_n(\tau_n)}}{c\tau_n}.
\end{equation*}
This implies that
\begin{equation}\label{st101}
\max_{1\le q\le m}\|\E\mathbf H^{q,c)}\|_2\le 2\max_{1\le q\le m}\|\E\mathbf X^{(q,c)}\|_2\le \frac{C\sqrt {L_n(\tau_n)}}{c\tau_n}.
\end{equation}
The inequalities (\ref{st100}) and (\ref{st101}) together imply that
\begin{equation}\label{stieltjes2}
 |s_n^{(c)}(\alpha,z)-{\ s}_n^{(c)}(\alpha,z)|\le  \frac{C\sqrt {L_n(\tau_n)}}{\sqrt n\tau_n v^2}.
\end{equation}
\end{proof}

According to Lemma \ref{truncation}  the matrices $\mathbf W$ and $\overline {\mathbf W}^{(c)}$ have the same limit distribution. In the what follows we shall assume without loss of generality that for any $n\ge 1$ and $q=1,\ldots, m$  and $j,k=1,\ldots,n$,
\begin{equation}\label{conditions}
 \E X^{(q)}_{jk}=0 \quad\text{and}\quad |X^{(q)}_{jk}|\le c\tau_n\sqrt n
\end{equation}
with $\tau_n\to 0$ such that
\begin{equation}\notag
 \qquad {L_n(\tau_n)}\to 0\qquad \text{and}\qquad\tau_n\sqrt n\to \infty\qquad{\text{as}}\qquad n\to \infty.
\end{equation}
We also have that
\begin{align}\label{var_1}
&\frac{1}{n^2} \sum_{j,k=1}^n |\E (X_{jk}^{(q)})^2- 1| \le C L_n(\tau_n),\\
\label{var_2}
&\frac{1}{n^2} \sum_{j,k=1}^n |\E X_{jk}^{(q)} X_{kj}^{(q)}- \rho| \le C L_n(\tau_n).
\end{align}

%%%%%%%%%%%%%%%%%%%%%%%%%%%%%%%%%%%%%%%%%%%%%%%%%%%%%%%%%%%%%%%%%%%%%%%%%%%%%%%%%%%%%%%%%%

\subsubsection{The universality of the limit distribution of singular values of shifted matrices}
In this Section we show that the limit distribution of singular values of product of random matrices satisfying assumptions of Theorem~\ref{shift} doesn't depend on the distribution of matrix entries. Let $\mathbf Y^{(1)}, \ldots,\mathbf Y^{(m)}$ be $n\times n$ independent random matrices  with independent Gaussian entries $n^{-1/2} Y_{jk}^{(q)}$ such that
\begin{align}
&\E Y_{jk}^{(q)}=0,\qquad \E (Y_{jk}^{(q)})^2=1,\qquad\text{for any}\quad q=1,\ldots,m,\, j,k=1\ldots,n;\notag\\
&\E Y_{jk}^{(q)} Y_{kj}^{(q)}=\rho\qquad\text{for any}\quad q=1,\ldots,m, 1\le j<k\le n.\notag
\end{align}
Vectors $(Y^{(q)}_{jk},Y^{(q)}_{kj})$ and r.v.'s $Y_{ll}^{(q)}$ for $q=1,\ldots,m$, $1\le j<k\le n$ and $l=1,\ldots,n$,
are mutually independent.
For any $\varphi\in[0,\frac{\pi}2]$ and any $\nu=1,\ldots,m$, introduce the matrices
\begin{equation}\notag
\mathbf Z^{(\nu)}(\varphi)=\mathbf X^{(\nu)}\cos\varphi+
\mathbf Y^{(\nu)}\sin\varphi
\end{equation}
where
$$[\mathbf Z^{(q)}(\varphi)]_{jk}=\frac1{\sqrt{n}}Z_{jk}^{(q)}=
\frac1{\sqrt{n}}(X_{jk}^{(q)}\cos\varphi+Y_{jk}^{(q}\sin\varphi).
$$
We define the matrices $\mathbf W(\varphi)$, $\mathbf H^{(q)}(\varphi)$, $\mathbf V(\varphi)$, $\widehat{\mathbf V}(\varphi)$, $\mathbf R(\varphi)$ by
\begin{align*}
\mathbf W(\varphi)=\prod_{\nu=1}^m\mathbf Z^{(\nu)}(\varphi),\quad
\mathbf H^{(\nu)}(\varphi)=\begin{bmatrix}&\mathbf Z^{(\nu)}(\varphi))&\mathbf  O\\&\mathbf O&\mathbf Z^{(m-\nu+1)}
(\varphi)\end{bmatrix}\\
\mathbf V(\varphi)=\prod_{\nu=1}^m\mathbf H^{(\nu)}(\varphi),\quad
\widehat{\mathbf V}(\varphi)=\mathbf V(\varphi)\mathbf J,\quad
\mathbf R(\varphi)=(\widehat{\mathbf V}(\varphi)-\mathbf J(z)-\alpha\mathbf I)^{-1}.
\end{align*}
Recall that $\mathbf I$ (with sub-index or without it) denotes the unit matrix of corresponding order, $\mathbf J(z)=\begin{bmatrix}&\mathbf O&
z\mathbf I\\ &\overline z\mathbf I&\mathbf O\end{bmatrix}$ and
$\mathbf O$ denotes the matrix with zero-entries.

In these notation the matrices $\mathbf W(0)$, $\mathbf H^{(\nu)}(0)$, $\mathbf V(0)$, $\widehat{\mathbf V}(0)$,
$\mathbf R(0)$ are generated by the matrices $\mathbf X^{(\nu)}$, $\nu=1,\ldots,m$, and
$\mathbf W(\frac{\pi}2)$, $\mathbf H^{(\nu)}(\frac{\pi}2)$, $\mathbf V(\frac{\pi}2)$,
$\widehat{\mathbf V}(\frac{\pi}2)$, $\mathbf R(\frac{\pi}2)$ are generated by $\mathbf Y^{(\nu)}$, $\nu=1,\ldots,m$. Let $s_n(\alpha,z,\varphi)$ denote the Stieltjes transform of symmetrized expected distribution function of singular values of the matrix $\mathbf W(\varphi)-z\mathbf I$. Then $s_n(\alpha,z,0)=s_n(\alpha,z)$ denote the Stieltjes transform of distribution function $G_n(x,z)$ and $s_n(\alpha,z,\frac{\pi}2)$ denote the Stieltjes transform of symmetrized expected distribution function of singular values of the matrix $\mathbf W(\frac{\pi}2)-z\mathbf I$ generated by $\mathbf Y^{(q)}$, $q=1,\ldots,m$. We prove the following Lemma.
\begin{lemma}\label{universality}
Under the assumptions of Theorem \ref{th:main} the following holds: for any $\delta>0$
\begin{equation*}
|s_n(\alpha,z,\frac{\pi}2)-s_n(\alpha,z,0)|\to 0\quad\text{as}\quad n\to\infty
\end{equation*}
uniformly in $\alpha=u+i v$ with $v\ge \delta$.
\end{lemma}
\begin{proof}
By Newton--Leibnitz formula we have
\begin{equation*}
s_n(\alpha,z,\frac{\pi}2)-s_n(\alpha,z,0)=\int_0^{\frac{\pi}2}\frac{\partial s_n(\alpha,z,\varphi)}{\partial \varphi}d\varphi.
\end{equation*}
Applying the formula for the derivative of matrix resolvent we get
\begin{equation}\label{d1}
\frac{\partial s_n(\alpha,z,\varphi)}{\partial \varphi}=-\frac1{2n}\E\Tr \mathbf R(\varphi)\frac{\partial \mathbf V(\varphi)}{\partial \varphi}
\mathbf J\mathbf R(\varphi).
\end{equation}
We shall omit in what follows the argument $\varphi$ in the notations of $\mathbf R$ and $\mathbf V$ if it doesn't confuse. By the definition of the matrix $\mathbf V$ and $\mathbf V_{a,b}$ (see~\eqref{vab}), we have
\begin{equation*}
 \frac{\partial \mathbf V}{\partial \varphi}=\sum_{q=1}^m\mathbf V_{1,q-1}
\frac{\partial\,\mathbf H^{(q)}}{\partial \varphi}\mathbf V_{q+1,m}.
\end{equation*}
Furthermore, by the definition of $\mathbf H^{(q)}$, for $q=1,\ldots, m$, we have
\begin{align}
\frac{\partial\,\mathbf H^{(q)}}{\partial \varphi}=\sum_{j=1}^{n}
\sum_{k=1}^{n}\Big(\frac{\partial\,\mathbf H^{(q)}}{\partial  Z^{(q)}_{jk}}\frac{d Z^{(q)}_{jk}}{d\varphi}
+\frac{\partial\,\mathbf H^{(q)}}
{\partial  Z^{(m-q+1)}_{jk}}\frac{d Z^{(m-q+1)}_{jk}}{d\varphi}
\Big),\notag
\end{align}
where we denote by $\mathbf e_j=(0,\ldots,0,1,\ldots,0)^T$ the column vector of the dimension $2n$ with all zero entries except $j$-th one, which equal to $1$,  $j=1,\ldots, 2n$. In these notations we have
\begin{align}
\frac{\partial\,\mathbf H^{(q)}}{\partial  Z^{(q)}_{jk}}=\frac1{\sqrt{n}}\mathbf e_j\mathbf e_k^T,\quad
\frac{\partial\,\mathbf H^{(q)}}{\partial  Z^{(m-q+1)}_{jk}}=\frac1{\sqrt{n}}\mathbf e_{k+n}\mathbf e_{j+n}^T,\notag
\end{align}
 for $j,k=1,\ldots, n$.
By the definition of $Z^{(q)}_{jk}$, we have
\begin{align}
\frac{dZ^{(q)}_{jk}}{d\varphi}=-X^{(q)}_{jk}\sin\varphi+ Y^{(q)}_{jk}\cos\varphi.\notag
\end{align}
After a simple calculation we get
\begin{align*}
\frac{\partial \mathbf V}{\partial \varphi}=\frac1{\sqrt{n}}\sum_{q=1}^m\sum_{j=1}^{n}
\sum_{k=1}^{n}\Bigg(\mathbf V_{1,q-1}\mathbf e_j{\mathbf e_k}^T\mathbf V_{q+1,m}(- X_{jk}^{(q)}\sin\varphi+ Y^{(q)}_{jk}\cos\varphi)\\
+
\mathbf V_{1,q-1}\mathbf e_{k+n}\mathbf e_{j+n}^T\mathbf V_{q+1,m}(- X_{jk}^{(m-q+1)}\sin\varphi+ Y^{(m-q+1)}_{jk}\cos\varphi)\Bigg).
\end{align*}
Introduce the following functions
\begin{align*}
u^{(q)}_{jk}&=-\Tr\mathbf R
\mathbf V_{1,q-1}\mathbf e_j\mathbf e_k^T\mathbf V_{q+1,m}\mathbf J\mathbf R,\\
v_{jk}^{(q)}&=\Tr\mathbf R\mathbf V_{1,q-1}\mathbf e_{k+n}{\mathbf e_{j+n}}^T\mathbf V_{q+1,m}\mathbf J\mathbf R,
\end{align*}
for $q=1,\ldots,m$, and  $j,k=1,\ldots,n$. In these notations we have
\begin{equation}\notag
\frac{\partial s_n(z,\varphi)}{\partial \varphi}=\Xi_1+\Xi_2,
\end{equation}
where
\begin{align*}
\Xi_1&=\frac1{2n\sqrt{n}}\sum_{q=1}^m\sum_{j=1}^{n}\sum_{k=1}^{n}
\E(-X_{jk}^{(q)}\sin\varphi+ Y^{(q)}_{jk}\cos\varphi))u^{(q)}_{jk}\\
\Xi_2=&\frac1{2n\sqrt{n}}\sum_{q=1}^m\sum_{j=1}^{n}\sum_{k=1}^{n}
\E(-X_{jk}^{(m-q+1)}\sin\varphi+ Y^{(m-q+1)}_{jk}\cos\varphi))v^{(q)}_{jk}.
\end{align*}
First we investigate $\Xi_1$. Let $\xi_{jk}^{(q)}= X_{jk}^{(q)}\cos\varphi+ Y_{jk}^{(q)}\sin\varphi$.
In  what follows we shall consider the functions
$u^{(q)}_{jk}= u^{(q)}_{jk}(\xi^{(q)}_{jk},\xi_{kj}^{(q)})$ as functions of $X_{jk}^{(q)},  X_{kj}^{(q)},
 Y_{jk}^{(q)}$ and $Y_{kj}^{(q)}$. Applying Taylor's formula, we may write
\begin{align}\label{te1}
u^{(q)}_{jk}(\xi^{(q)}_{jk},\xi_{kj}^{(q)})&=u^{(q)}_{jk}(0,0)+\xi_{jk}^{(q)}\frac{\partial u^{(q)}_{jk}}{\partial \xi_{jk}^{(q)}}(0,0)
+\xi_{kj}^{(q)}\frac{\partial u^{(q)}_{jk}}{\partial \xi_{kj}^{(q)}}(0,0)
\notag\\
&+\E_{\theta}(\xi_{jk}^{(q)})^2
(1-\theta)\frac{\partial^2u_{jk}^{(q))}}{\partial{\xi_{jk}^{(q)}}^2}(\theta \xi_{jk}^{(q)}, \theta \xi_{kj}^{(q)})\notag\\&+
2\E_{\theta}\xi_{kj}^{(q)}\xi_{jk}^{(q)}\frac{\partial^2 u_{jk}^{(q)}}{\partial \xi_{jk}^{(q)}\partial \xi_{kj}^{(q)}}(\theta \xi_{jk}^{(q)},\theta\xi_{kj}^{(q)})
\notag\\&+\E_{\theta}(\xi_{kj}^{(q)})^2(1-\theta)\frac{\partial^2u_{jk}^{(q))}}{\partial{\xi_{kj}^{(q)}}^2}(\theta \xi_{jk}^{(q)}, \theta \xi_{kj}^{(q)}).
\end{align}
Here $\theta$ are uniformly distributed on $[0,1]$ and is independent of all $X_{jk}^{(q)}$ and $Y_{jk}^{(q)}$, and $\E_{\theta}$ denotes the expectation with respect to the random variable $\theta$. Furthermore, we introduce the random variables
\begin{align}
{{\widehat\xi_{jk}}^{(q)}}&=- X_{jk}^{(q)}\sin\varphi+ Y^{(q)}_{jk}\cos\varphi.\notag
\end{align}
Multiplying~\eqref{te1} by $\widehat \xi_{jk}^{(q)}$  and taking expectation, we rewrite $\Xi_1$ as $\Xi_1 = \Xi_{11} + \Xi_{12}$, where
\begin{align*}
\Xi_{11} &=\E {{\widehat\xi_{jk}}^{(q)}}\xi_{jk}^{(q)}
\E \frac{\partial u^{(q)}_{jk}}{\partial \xi_{jk}^{(q)}}(0,0)
+\E{{\widehat\xi_{jk}}^{(q)}}\xi_{kj}^{(q)}\E\frac{\partial u^{(q)}_{jk}}{\partial \xi_{kj}^{(q)}}(0,0),\\
\Xi_{12} &= \E{{\widehat\xi_{jk}}^{(q)}}(\xi_{jk}^{(q)})^2
(1-\theta)\frac{\partial^2u_{jk}^{(q))}}{\partial{{{\xi_{jk}^{(q)}}}}^2}(\theta {{\xi_{jk}^{(q)}}},\theta \xi_{kj}^{(q)})\\
&+ 2\E{{\widehat\xi_{jk}}^{(q)}}\xi_{jk}^{(q)}\xi_{kj}^{(q)}\frac{\partial^2 u_{jk}^{(q)}}{\partial \xi_{jk}^{(q)}\partial
\xi_{kj}^{(q)}}(\theta\xi_{jk}^{(q)},\theta\xi_{kj}^{(q)})\\
&+\E(\xi_{kj}^{(q)})^2{{\widehat\xi_{jk}}^{(q)}}(1-\theta)
\frac{\partial^2u_{jk}^{(q))}}{\partial{\xi_{kj}^{(q)}}^2}(\theta \xi_{jk}^{(q)}, \theta \xi_{kj}^{(q)}).
\end{align*}
It is straightforward to check, that
\begin{align*}
&\E {{\widehat\xi_{jk}}^{(q)}}\xi_{jk}^{(q)}= \cos \varphi \sin \varphi \E[(Y_{jk}^{(q)})^2 - (X_{jk}^{(q)})^2]\\
&\E{{\widehat\xi_{jk}}^{(q)}}\xi_{kj}^{(q)}= \cos \varphi \sin \varphi \E[Y_{jk}^{(q)} Y_{kj}^{(q)} - X_{jk}^{(q)}X_{kj}^{(q)}]
\end{align*}
We introduce the following matrices
\begin{align}
\mathbf B_{jk}^{(q)}:=&\mathbf V_{1,q-1}\mathbf e_j{\mathbf e_k}^T\mathbf V_{q+1,m}.\notag
\end{align}
In these notations we get $u^{(q)}_{jk}=-\Tr\mathbf B_{jk}^{(q)}\mathbf J\mathbf R^2$. It is easy to check that
\begin{align*}
\frac{\partial u_{jk}^{(q)}}{\partial{\xi_{jk}^{(q)}}}(\theta_1\xi_{jk}^{(q)},\theta\xi_{kj}^{(q)})&=
-\Tr \frac{\partial \mathbf B_{jk}^{(q)}}{\partial{\xi_{jk}^{(q)}}}\mathbf J\mathbf R^2
+\Tr \mathbf B_{jk}^{(q)}\mathbf J\mathbf R^2
\frac{\partial \mathbf V}{\partial \xi_{jk}^{(q)}}\mathbf J\mathbf R\\
&+\Tr\mathbf B_{jk}^{(q)}\mathbf J\mathbf R\frac{\partial \mathbf V}{\partial \xi_{jk}^{(q)}}\mathbf J
\mathbf R^2 = I_1+I_2+I_3.
\end{align*}
Furthermore,
\begin{align*}
\frac{\partial \mathbf B_{jk}^{(q)}}{\partial \xi_{jk}^{(q)}}&=
\frac1{\sqrt{n}}\mathbf V_{1,m-q}\mathbf e_{k+n}\mathbf e_{j+n}^T
\mathbf V_{m-q+2,q-1}\mathbf e_j\mathbf e_k^T\mathbf V_{q+1,m}\mathbb I\{m-q\le q-1\}\\
&+\frac1{\sqrt{n}}\mathbf V_{1,q-1}
\mathbf e_j{\mathbf e_k}^T\mathbf V_{q+1,m-q}\mathbf e_{k+n}
\mathbf e_{j+n}^T\mathbf V_{m-q+2,m}\mathbb I\{m-q\ge  q\},
\end{align*}
and
\begin{equation}\label{der4}
\frac{\partial \mathbf V}{\partial \xi_{jk}^{(q)}}=
\frac1{\sqrt{n}}\mathbf V_{1,q-1}\mathbf e_j\mathbf e_k^T\mathbf V_{q+1,m}+\frac1{\sqrt{n}}\mathbf V_{1,m-q}
\mathbf e_{k+n}\mathbf e_{j+n}^T\mathbf V_{m-q+2,m}.
\end{equation}

%%%%%%%%%%%%%%%%%%%%%%%%%%%%%%%%%%%%%%%%%%%%%%%%%%%%%%%%%%%%%%%%%%%
Note that $[\mathbf V_{m-q+2, q-1}]_{j,j+n}=0$ and
$[\mathbf V_{q+1,m-q}]_{k+n,k}=0$. These equalities imply that
\begin{align}\notag
I_{1}=0.
\end{align}
Using~\eqref{der4} we get
\begin{align}\label{T212}
I_2=I_{21}+I_{22},
\end{align}
where
\begin{align*}
I_{21}&=\frac1{\sqrt n}\ \Tr\mathbf V_{1,q-1}\mathbf e_j{\mathbf e_k}^T\mathbf V_{q+1,m}\mathbf J\mathbf R^2 \mathbf V_{1,q-1}\mathbf e_j{\mathbf e_k}^T\mathbf V_{q+1,m}\mathbf J
\mathbf R,\notag\\
I_{22}&=\frac1{\sqrt n}\ \Tr\mathbf V_{1,q-1}\mathbf e_j{\mathbf e_k}^T\mathbf V_{q+1,m}\mathbf J\mathbf R^2 \mathbf V_{1,m-q}
\mathbf e_{k+n}{\mathbf e_{j+n}}^T\mathbf V_{m-q+2,m}\mathbf J\mathbf R
\end{align*}
We shall bound each term in~\eqref{T212}. Note that
\begin{align}
I_{21}=\frac1{\sqrt n}[\mathbf V_{q+1,m}\mathbf J\mathbf R^2\mathbf V_{1,q-1}]_{kj}
[\mathbf V_{q+1,m}\mathbf J\mathbf R\mathbf V_{1,q-1}]_{kj}.\notag
\end{align}
It is straightforward to check that
\begin{align}
|I_{21}|\le C\,v^{-3}n^{-1/2}\|{\mathbf e_{k}}^T\mathbf V_{q+1,m}\|_2^2 \|\mathbf V_{1,q-1}{\mathbf e_{j}}\|_2^2. \notag
\end{align}
Note that the random variables in the r.h.s of the last inequality conditionally independent with respect to $\xi_{jk}^{(q)}$ and $\xi_{kj}^{(q)}$. We may write
\begin{align*}
\E\Big\{| I_{21}|\Big|\xi_{jk}^{(q)},\xi_{kj}^{(q)}\Big\}&\le \\
&\frac{C}{v^{3}n^{1/2}}\E\Big\{\|{\mathbf e_{k}}^T\mathbf V_{q+1,m}\|_2^2 \Big|\xi_{jk}^{(q)},\xi_{kj}^{(q)}\Big\}\E\Big\{\|\mathbf V_{1,q-1}{\mathbf e_{j}}\|_2^2\Big|\xi_{jk}^{(q)},\xi_{kj}^{(q)}\Big\}.
\end{align*}
Applying Lemma \ref{norm4}, we get
\begin{equation*}
\E\Big\{|I_{21}|\Big|\xi_{jk}^{(q)},\xi_{kj}^{(q)}\Big\}\le Cn^{-1/2}v^{-3}.
\end{equation*}
Similarly we estimate $I_{22}$ and $I_3$. It follows from these bounds,~\eqref{var_1} and~\eqref{var_2} that
\begin{equation*}
|\Xi_{11}|\le C v^{-3} L_n(\tau_n).
\end{equation*}

We now estimate $\Xi_{12}$. Without loss of generality we may assume that
\begin{equation*}
\max\Big\{|\xi_{jk}^{(q)}|,|\xi_{kj}^{(q)}|,|{\widehat{\xi}}_{jk}^{(q)}|,
|{\widehat{\xi}}_{kj}^{(q)}|\Big\}\le C\tau_n\sqrt n.
\end{equation*}
If we prove that there exists a constant $C$ such that, for any $q=1,\ldots,m$, $1\le j,k\le n$,
\begin{align}\label{w1}
\max\Big\{\Big|\E\Big\{\frac{\partial^2 u_{jk}^{(q)}}{\partial{\xi_{jk}^{(q)}}^2}(\theta\xi_{jk}^{(q)},\theta\xi_{kj}^{(q)})
\Big|\xi_{jk}^{(q},\xi_{kj}^{(q)}\Big\}\Big|,\Big|\E\{\frac{\partial^2 u_{jk}^{(q)}}{\partial \xi_{jk}^{(q)}
\partial \xi_{kj}^{(q)}}(\theta\xi_{jk}^{(q)},\theta\xi_{kj}^{(q)})\Big|\xi_{jk}^{(q},\xi_{kj}^{(q)}\Big\}\Big|\Big\}\nonumber\\
\le Cn^{-1}v^{-4},
\end{align}
we get
\begin{equation*}
|\E {{\widehat\xi_{jk}}^{(q)}}
u^{(q)}_{jk}(\xi_{jk}^{(q)},\xi_{kj}^{(q)})|\le \frac{C\tau_n}{\sqrt n}.
\end{equation*}
The last bound implies that
\begin{equation}\label{xi1}
|\Xi_{12}|\le C\tau_nv^{-4}.
\end{equation}

Furthermore,
\begin{align*}
%\label{der1}
&\frac{\partial^2 u_{jk}^{(q)}}{\partial{\xi_{jk}^{(q)}}^2}(\theta \xi_{jk}^{(q)},\theta\xi_{kj}^{(q)})=
-2\Tr\mathbf B_{jk}^{(q)}\mathbf J\RR^2\frac{\partial \V}{\partial \xi_{jk}^{(q)}}\mathbf J
\RR\frac{\partial \V}{\partial \xi_{jk}^{(q)}}\mathbf J\RR \\
&-2\Tr\mathbf B_{jk}^{(q)}\mathbf J\RR\frac{\partial \V}{\partial \xi_{jk}^{(q)}}\mathbf J
\RR^2\frac{\partial \V}{\partial \xi_{jk}^{(q)}}\mathbf J\RR-2\Tr\mathbf B_{jk}^{(q)}\mathbf J\RR\frac{\partial \V}{\partial \xi_{jk}^{(q)}}\mathbf J
\RR\frac{\partial \V}{\partial \xi_{jk}^{(q)}}\mathbf J\RR^2\\
&=T_1+T_2+T_3.
\end{align*}

We bound $T_1$ now. The estimates for $T_2, T_3$ may be written down in the similar way. Using~\eqref{der4} we get
\begin{align}\label{T1234}
T_1=T_{11}+\cdots+T_{14},
\end{align}
where
\begin{align*}
T_{11}&=-2\frac1{n}\ \Tr\mathbf V_{1,q-1}\mathbf e_j{\mathbf e_k}^T\mathbf V_{q+1,m}\mathbf J\mathbf R^2\\
&\qquad\qquad\times\mathbf V_{1,q-1}\mathbf e_j{\mathbf e_k}^T\mathbf V_{q+1,m}\mathbf J
\mathbf R\mathbf V_{1,q-1}\mathbf e_j{\mathbf e_k}^T\mathbf V_{q+1,m}\mathbf J\mathbf R,\\
T_{12}&=-2\frac1{n}\ \Tr\mathbf V_{1,q-1}\mathbf e_j{\mathbf e_k}^T\mathbf V_{q+1,m}\mathbf J\mathbf R^2\\
&\qquad\qquad\times\mathbf V_{1,q-1}\mathbf e_j{\mathbf e_k}^T\mathbf V_{q+1,m}\mathbf J\mathbf R\mathbf V_{1,m-q}
\mathbf e_{k+n}{\mathbf e_{j+n}}^T\mathbf V_{m-q+2,m}\mathbf J\mathbf R,\notag\\
T_{13}&=-2\frac1{n}\ \Tr\mathbf V_{1,q-1}\mathbf e_j{\mathbf e_k}^T\mathbf V_{q+1,m}\mathbf J\mathbf R^2\\
&\qquad\qquad\times\mathbf \mathbf V_{1,m-q}
\mathbf e_{k+n}{\mathbf e_{j+n}}^T\mathbf V_{m-q+2,m}\mathbf J\mathbf R
\mathbf V_{1,q-1}\mathbf e_j{\mathbf e_k}^T\mathbf V_{q+1,m}\mathbf J\mathbf R,\\
T_{14}&=-2\frac1{n}\ \Tr\mathbf V_{1,q-1}\mathbf e_j{\mathbf e_k}^T\mathbf V_{q+1,m}\mathbf J\mathbf R^2\\
&\qquad\qquad\times \mathbf V_{1,m-q}
\mathbf e_{k+n}{\mathbf e_{j+n}}^T\mathbf V_{m-q+2,m}\mathbf J
\mathbf R\mathbf V_{1,m-q}
\mathbf e_{k+n}{\mathbf e_{j+n}}^T\mathbf V_{m-q+2,m}\mathbf J\mathbf R,
\end{align*}
We shall bound each term in~\eqref{T1234}. Note that
\begin{align}
T_{11}=-2\frac1{n}[\mathbf V_{q+1,m}\mathbf J\mathbf R^2\mathbf V_{1,q-1}]_{kj}
[\mathbf V_{q+1,m}\mathbf J\mathbf R\mathbf V_{1,q-1}]_{kj}[\mathbf V_{q+1,m}\mathbf J\mathbf R\mathbf V_{1,q-1}]_{kj}.\notag
\end{align}
It is straightforward to check that
\begin{align}
|T_{31}|\le C\,v^{-4}n^{-1}\|{\mathbf e_{k}}^T\mathbf V_{q+1,m}\|_2^3 \|\mathbf V_{1,q-1}{\mathbf e_{j}}\|_2^3. \notag
\end{align}
Note that the random variables in the r.h.s of the last inequality conditionally independent with respect to $\xi_{jk}^{(q)}$ and $\xi_{kj}^{(q)}$. We may write
\begin{align}\label{T31}
\E\Big\{| T_{31}|\Big|\xi_{jk}^{(q)},\xi_{kj}^{(q)}\Big\}&\le \nonumber \\
&\frac{C}{v^{4}n}\E\Big\{\|{\mathbf e_{k}}^T\mathbf V_{q+1,m}\|_2^3
\Big|\xi_{jk}^{(q)},\xi_{kj}^{(q)}\Big\}\E\Big\{\|\mathbf V_{1,q-1}{\mathbf e_{j}}\|_2^3\Big|\xi_{jk}^{(q)},\xi_{kj}^{(q)}\Big\}.
\end{align}
Applying Lemma \ref{norm4}, we get
\begin{equation*}
\E\Big\{|T_{31}|\Big|\xi_{jk}^{(q)},\xi_{kj}^{(q)}\Big\}\le Cn^{-1}v^{-4}.
\end{equation*}
Furthermore, we represent $T_{32}$ in the form
\begin{align*}
T_{32}=-2\frac1{n}[\mathbf V_{q+1,m}\mathbf J\mathbf R^2\mathbf V_{1,q-1}]_{k,j}
[\mathbf V_{q+1,m}\mathbf J\mathbf R\mathbf V_{1,m-q}]_{k,k+n}
[\mathbf V_{m-q+2,m}\mathbf J\mathbf R\mathbf V_{1,q-1}]_{j+n,j}.
\end{align*}
Similar to (\ref{T31}) we get
\begin{align*}
|T_{32}|\le C\,v^{-4}n^{-1}\,\|{\mathbf e_k}^T\mathbf V_{q+1,m}\|_2^2\,
\|\mathbf V_{1,q-1}\mathbf e_j\|_2^2\\ \times
\|\mathbf V_{1,m-q}\mathbf e_{k+n}\|_2\,\|{\mathbf e_{j+n}}^T\mathbf V_{m-q+2,m}\|_2.
\end{align*}
Applying H\"older's inequality, we get
\begin{align*}
\E\Big\{| T_{32}|\Big|\xi_{jk}^{(q)},\xi_{kj}^{(q)}\Big\}\le \frac{C}{v^{4}n} \E^{\frac12}\Big\{\|{\mathbf e_k}^T
\mathbf V_{q+1,m}\|_2^4\Big|\xi_{jk}^{(q)},\xi_{kj}^{(q)}\Big\}\E^{\frac12}\Big\{\|\mathbf V_{1,q-1}\mathbf e_j\|_2^4
\Big|\xi_{jk}^{(q)},\xi_{kj}^{(q)}\Big\}\\
\times\E^{\frac12}\Big\{\|\mathbf V_{1,m-q}\mathbf e_{k+n}\|_2^2\Big|\xi_{jk}^{(q)},\xi_{kj}^{(q)}\Big\}
\E^{\frac12}\Big\{\|{\mathbf e_{j+n}}^T\mathbf V_{m-q+2,m}\|_2^2\Big|\xi_{jk}^{(q)},\xi_{kj}^{(q)}\Big\}.
\end{align*}
Using Lemma \ref{norm4}, we get
\begin{equation*}
\E\Big\{| T_{12}|\Big|\xi_{jk}^{(q)},\xi_{kj}^{(q)}\Big\}\le C\,v^{-4}n^{-1}.
\end{equation*}
Analogously we get the bounds for other terms $T_{1l}$, for $l=3,4$. We have
\begin{equation*}
\E\Big\{| T_{1}|\Big|\xi_{jk}^{(q)},\xi_{kj}^{(q)}\Big\}\le C\,v^{-4}n^{-1}.
\end{equation*}
This proves~\eqref{w1} and~\eqref{xi1}. Similarly we may estimate the term $\Xi_2$
\begin{equation*}
|\Xi_2|\le C\tau_nv^{-4}.
\end{equation*}
It follows that there exists some $\delta > 0$ such that
\begin{equation*}
\lim_{n \rightarrow \infty}|s_n(\alpha,z,\frac{\pi}2)-s_n(\alpha,z,0)| = 0,
\end{equation*}
for all $v \geq \delta$.
The last inequality proves the Lemma~\ref{universality}.
\end{proof}

\subsubsection{ The Limit Distribution of Singular Values of $\mathbf V(z)$ in the Gaussian case} {In this Section we find the limit distribution of singular values of shifted products of Gaussian random matrices.}
Recall that
\begin{align*}
\mathbf H^{(\nu)}=\left(\begin{matrix}&\mathbf Y^{(\nu)}& \mathbf O\\
&\mathbf O & \mathbf {Y^{(m-\nu+1)}}^*\end{matrix}\right), \quad
\mathbf J(z):=\left(\begin{matrix}&\mathbf O  & z\ \mathbf I \\&\overline z\ \mathbf I &\mathbf O\end{matrix}\right),
\text{ and } \mathbf J:=\mathbf J(1).
\end{align*}
For any $1\le a, b\le m$, put
$$
\mathbf V_{[a,b]}=\begin{cases}\prod_{k=a}^{b}\mathbf H^{(k)},\quad\text{for }a\le b,\\
                   \mathbf I\quad\text{otherwise},
                  \end{cases}
$$
and
$$
\mathbf V(z):=\mathbf V\mathbf J-\mathbf J(z),\quad \mathbf R=(\mathbf V(z)-\alpha \mathbf I)^{-1}.
$$
It is straightforward to check
\begin{align}\label{systemdef}
s_n(\alpha,z)&=\frac1n\sum_{j=1}^n\E[\mathbf R(\alpha,z)]_{jj}\nonumber\\
&=\frac1n\sum_{j=1}^n\E[\mathbf R(\alpha,z)]_{j+nj+n}=\frac1{2n}\sum_{j=1}^{2n}\E[\mathbf R(\alpha,z)]_{jj}.
\end{align}
We introduce the following functions
\begin{align*}
t_n(\alpha,z)&=\frac1n\sum_{j=1}^n\E[\mathbf R(\alpha,z)]_{j+n,j},\quad
u_n(\alpha,z)=\frac1n\sum_{j=1}^n\E[\mathbf R(\alpha,z)]_{j,j+n}.
\end{align*}
We prove the following statement
\begin{statement}\label{respect}
Let r.v.'s  $Y_{jk}^{(q)}$, $q=1,\ldots,m$, $j,k=1,\ldots n$ are Gaussian and satisfy the conditions $\Cond$.
Then the following  limit exists
$$
g=g(\alpha,z)=\lim_{n\to\infty}s_n(\alpha,z),
$$
and  satisfy the system equations
\begin{align}\label{system00}
 &1+wg+(-1)^{m+1}w^{m-1}g^{m+1}=0,\nonumber\\
&g(w-\alpha)^2+(w-\alpha)-g|z|^2=0,
\end{align}
{with a function $w=w(\alpha,z)$ such that $\imag(w-\alpha)>0$}.
\end{statement}
\begin{corollary}\label{cor}Under the assumptions of Theorem \ref{shift} for any $z\in\mathbb C$  there exists a distribution function $G(x,z)$ such that
$\lim_{n\to\infty} G_n(x,z)=G(x,z)$ and ${g}={g}(\alpha,z)=\int_{-\infty}^{\infty}\frac1{x-\alpha}dG(x,z)$
satisfy the system of equations~\eqref{system00} and
\begin{equation}\label{distance}
\Delta_n(z):=\sup_x|G_n(x,z)-G(x,z)|\to 0\quad\text{as}\quad n\to \infty.
\end{equation}
\end{corollary}
\begin{remark}\label{rem31}
Note that the second equation of~\eqref{system} implies
\begin{equation}\notag
\imag{g}=-\imag\Big\{\frac{w-\alpha}{(w-\alpha)^2-|z|^2}\Big\}=\frac{\imag\{w-\alpha\}(|w-\alpha|^2+|z|^2)}{|(w-\alpha)^2-|z|^2|^2}.
\end{equation}
This equality implies that $\imag(w-\alpha)>0$.
\end{remark}

\begin{proof} {\it Statement} \ref{respect}. In what follows we shall denote by $\varepsilon_n(\alpha,z)$
 a generic  error function such that
$|\varepsilon_n(\alpha,z)|\le \frac{C\tau_n^q}{v^r}$ for some positive constants $C,q,r$. By the resolvent equality, we may write
\begin{align}\label{b1}
 1+\alpha s_n(\alpha,z)&=\frac1{2n}\E\Tr\mathbf V(z)\mathbf R(\alpha,z) \nonumber \\
&= \frac1{2n}\E\Tr \mathbf V\mathbf J\mathbf R(\alpha,z)-\frac12zt_n(\alpha,z)-\frac12\overline z u_n(\alpha,z).
\end{align}
In the following we shall write $\mathbf R$ instead of $\mathbf R(\alpha,z)$.
Introduce the notation
\begin{equation}\notag
 \mathcal A:=\frac1{2n}\E\Tr \mathbf V\mathbf J\mathbf R
\end{equation}
and {represent $\mathcal A$ as follows}
\begin{equation}\notag
 \mathcal A=\frac12\mathcal A_1+\frac12\mathcal A_2,
\end{equation}
where
$$
\mathcal A_1=\frac1n\sum_{j=1}^n\E[\mathbf V\mathbf J\mathbf R]_{jj},\quad
\mathcal A_2=\frac1n\sum_{j=1}^n\E[\mathbf V\mathbf J\mathbf R]_{j+n,j+n}.
$$
By definition of the matrix $\mathbf V$ {and the matrix $\mathbf H^{(1)}$}, we have
\begin{equation}\notag
 \mathcal A_1=\frac1{n\sqrt n}\sum_{j,k=1}^n\E Y^{(1)}_{jk}[\mathbf V_{2,m}\mathbf J\mathbf R]_{kj}.
\end{equation}
In the Gaussian case we may represent the random variables $Y_{jk}^{(q)}$ and $Y_{kj}^{(q)}$ in the form
\begin{align}\label{gausscorel}
Y_{jk}^{(q)}&=a\xi_{jk}^{(q)}+b\eta_{jk}^{(q)},\notag\\
Y_{kj}^{(q)}&=a\xi_{jk}^{(q)}-b\eta_{jk}^{(q)},
\end{align}
where $a=\sqrt{\frac{1+\rho}2}$, $b=\sqrt{\frac{1-\rho}2}$ and $\xi_{jk}^{(q)}$, $\eta_{jk}^{(q)}$
are mutually independent standard Gaussian r.v.'s.
We shall use the well-known equality for the standard Gaussian r.v. $\xi$ and any smooth function $f$
\begin{equation}\label{gauss}
 \E\xi f(\xi)=\E f'(\xi).
\end{equation}
First we represent $\mathcal A_1$ in the form
\begin{equation*}
 \mathcal A_1=\mathcal A_{11}+\mathcal A_{12}+\mathcal A_{13},
\end{equation*}
where
\begin{align*}
 \mathcal A_{11}&=\frac1{n\sqrt n}\sum_{j=1}^{n-1}\sum_{k=j+1}^n\E Y_{jk}^{(1)}[\mathbf V_{2,m}\mathbf J\mathbf R]_{kj},\\
\mathcal A_{12}&=\frac1{n\sqrt n}\sum_{j=1}^{n}\E Y_{jj}^{(1)}[\mathbf V_{2,m}\mathbf J\mathbf R]_{jj},\\
\mathcal A_{13}&=\frac1{n\sqrt n}\sum_{j=2}^{n}\sum_{k=1}^{j-1}\E Y_{jk}^{(1)}[\mathbf V_{2,m}\mathbf J\mathbf R]_{kj}.
\end{align*}
First we note that
\begin{align*}
 |\mathcal A_{12}|&\le \frac1{n\sqrt n}\sum_{j=1}^n \E^{\frac12}|\mathbf e_j^T\mathbf V_{2,m}\mathbf J\mathbf R\mathbf e_j|^2\le
\frac1{\sqrt n}\Big(\frac1n\sum_{j=1}^n\E|\mathbf e_j^T\mathbf V_{2,m}\mathbf J\mathbf R\mathbf e_j|^2\Big)^{\frac12}\\
&\le \frac1{v\sqrt n}\Big(\frac1n\sum_{j=1}^n\E\|\mathbf e_j^T\mathbf V_{2,m}\|^2\Big)^{\frac12}\le \frac C{v\sqrt n}.
\end{align*}
We use here the inequalities $\|\mathbf J\mathbf R\mathbf e_j\|\le \|\mathbf J\mathbf R\|\le v^{-1}$ and
$ |\mathbf e_j^T\mathbf V_{2,m}\mathbf J\mathbf R\mathbf e_j|\le \||\mathbf e_j^T\mathbf V_{2,m}\|\|\mathbf J\mathbf R\mathbf e_j\|.$
We may write now
\begin{equation}\label{a12}
 \mathcal A_{12}=\varepsilon_n(\alpha,z).
\end{equation}
Furthermore, we consider $\mathcal A_{11}$ and $\mathcal A_{13}$. Using~\eqref{gausscorel}, we get
\begin{align*}
 \mathcal A_{11}=\frac1{n\sqrt n}\sum_{j=1}^{n-1}\sum_{k=j+1}^n(a\E \xi_{jk}^{(1)}[\mathbf V_{2,m}\mathbf J\mathbf R]_{kj}+b
\E \eta_{jk}^{(1)}[\mathbf V_{2,m}\mathbf J\mathbf R]_{kj}).
\end{align*}
Applying~\eqref{gauss}, we get
\begin{equation*}
 \mathcal A_{11}=\frac1{n\sqrt n}\sum_{j=1}^{n-1}\sum_{k=j+1}^n\left(a\E\left[\frac{\partial \mathbf V_{2,m}\mathbf J\mathbf R}
{\partial \xi_{jk}^{(1)}}\right]_{kj}+b\E\left[\frac{\partial \mathbf V_{2,m}\mathbf J\mathbf R}
{\partial \eta_{jk}^{(1)}}\right]_{kj}\right).
\end{equation*}
A simple calculation shows that
\begin{equation*}
 \mathcal A_{13}=\frac1{n\sqrt n}\sum_{j=1}^{n-1}\sum_{k=j+1}^n( a\E\xi_{jk}^{(1)}[\mathbf V_{2,m}\mathbf J\mathbf R]_{kj}
-b\E\eta_{jk}^{(1)}[\mathbf V_{2,m}\mathbf J\mathbf R]_{jk}).
\end{equation*}
By the equality (\ref{gauss}), we have
\begin{equation*}
 \mathcal A_{13}=\frac 1{n\sqrt n}\sum_{j=1}^{n-1}\sum_{k=j+1}^n\Big(a\E\Big[\frac{\partial \mathbf V_{2,m}\mathbf J\mathbf R}
{\partial \xi_{jk}^{(1)}}
\Big]_{jk}-b\E\Big[\frac{\partial \mathbf V_{2,m}\mathbf J\mathbf R}{\partial \eta_{jk}^{(1)}}
\Big]_{jk}\Big).
\end{equation*}

Note that for $1\le j<k\le n$
\begin{align}\label{deriv1}
 \frac{\partial \mathbf V_{2,m}\mathbf J\mathbf R}
{\partial \xi_{jk}^{(1)}}=a\left(\frac{\partial \mathbf V_{2,m}\mathbf J\mathbf R}
{\partial Y_{jk}^{(1)}}+\frac{\partial \mathbf V_{2,m}\mathbf J\mathbf R}
{\partial Y_{kj}^{(1)}}\right),\nonumber\\
\frac{\partial \mathbf V_{2,m}\mathbf J\mathbf R}
{\partial \eta_{jk}^{(1)}}=b\left(\frac{\partial \mathbf V_{2,m}\mathbf J\mathbf R}
{\partial Y_{jk}^{(1)}}-\frac{\partial \mathbf V_{2,m}\mathbf J\mathbf R}
{\partial Y_{kj}^{(1)}}\right).
\end{align}
Computing the matrix derivatives
\begin{align}\label{deriv2}
 \frac{\partial \mathbf V_{2,m}\mathbf J\mathbf R}{\partial Y^{(1)}_{jk}}&=\frac1{\sqrt n}
\mathbf V_{2,m-1}\mathbf e_{k+n}\mathbf e_{j+n}^T\mathbf J\mathbf R\nonumber\\
&-\frac1{\sqrt n}\mathbf V_{2,m}
\mathbf J\mathbf R\mathbf e_j\mathbf e_k^T\mathbf V_{2,m}\mathbf J\mathbf R-\mathbf V_{2,m}
\mathbf J\mathbf R\mathbf V_{1,m-1}\mathbf e_{k+n}\mathbf e_{j+n}^T\mathbf J\mathbf R,\nonumber\\
\frac{\partial \mathbf V_{2,m}\mathbf J\mathbf R}{\partial Y^{(1)}_{kj}}&=\frac1{\sqrt n}
\mathbf V_{2,m-1}\mathbf e_{j+n}\mathbf e_{k+n}^T\mathbf J\mathbf R\nonumber\\
&-\frac1{\sqrt n}\mathbf V_{2,m}
\mathbf J\mathbf R\mathbf e_{k}\mathbf e_{j}^T\mathbf V_{2,m}\mathbf J\mathbf R-\mathbf V_{2,m}
\mathbf J\mathbf R\mathbf V_{1,m-1}\mathbf e_{j+n}\mathbf e_{k+n}^T\mathbf J\mathbf R.
\end{align}
Combining the equalities~\eqref{deriv1} and~\eqref{deriv2}, we get
\begin{align*}
\frac{\partial \V_{2,m}\J\RR}{\partial \xi_{jk}^{(1)}}=&a\frac{1}{\sqrt n} ( \V_{2,m-1}(\ee_{k+n}\ee_{j+n}^T +\ee_{j+n}\ee_{k+n}^T)\J\RR \\
&-\V_{2,m}\J\RR(\ee_j\ee_k^T+\ee_k\ee_j^T)\V_{2,m}\J\RR \\ &-\V_{2,m}\J\RR\V_{1,m-1}(\ee_{k+n}\ee_{j+n}^T+\ee_{j+n}\ee_{k+n}^T)\J\RR)\\
\frac{\partial \V_{2,m}\J\RR}{\partial \eta_{jk}^{(1)}}=&-b\frac1{\sqrt n} (\V_{2,m-1}(\ee_{k+n}\ee_{j+n}^T-\ee_{j+n}\ee_{k+n}^T)\J\RR\\
&-\V_{2,m}\J\RR(\ee_j\ee_k^T-\ee_{k}\ee_{j}^T)\V_{2,m}\J\RR \\
&-\V_{2,m}\J\RR\V_{1,m-1}(\ee_{k+n}\ee_{j+n}^T-\ee_{j+n}\ee_{k+n}^T)\J\RR).
\end{align*}
Using the previous steps we may write
\begin{equation}\label{a11}
 \mathcal A_{11}+\mathcal A_{13}=\mathcal A_{111}+\ldots+\mathcal A_{114},
\end{equation}
where
\begin{align*}
 \mathcal A_{111}&=-\frac {2\rho}{n^2}\sum_{j=1}^{n-1}\sum_{k=j+1}^n\E[\V_{2,m}\J\RR]_{jj}[\V_{2,m}\J\RR]_{kk},\\
\mathcal A_{112}&=-\frac {1}{n^2}\sum_{j=1}^{n-1}\sum_{k=j+1}^n\E\Big([\V_{2,m}\J\RR]_{jk}^2
+[\V_{2,m}\J\RR]_{kj}^2\Big),\notag\\
\mathcal A_{113}&=-\frac {1}{n^2}\sum_{j=1}^{n-1}\sum_{k=j+1}^n\Big(\E[\V_{2,m}\J\RR \V_{1,m-1}]_{kk+n}[\J\RR]_{j+n,j} \\
&\qquad\qquad\qquad\qquad+\E[\V_{2,m}\J\RR\V_{1,m-1}]_{j,j+n}[\J\RR]_{k+n,k}\Big),\\
\mathcal A_{114}&=-\frac {\rho}{n^2}\sum_{j=1}^{n-1}\sum_{k=j+1}^n\Big(\E[\V_{2,m}\J\RR \V_{1,m-1}]_{k,j+n}[\J\RR]_{k+n,j}\\
&\qquad\qquad\qquad\qquad+\E[\V_{2,m}\J\RR\V_{1,m-1}]_{j,k+n}[\J\RR]_{j+n,k}\Big).
\end{align*}
We use here $[V_{2,m-1}]_{k,k+n}=[V_{2,m-1}]_{j,k+n}=[V_{2,m-1}]_{k,j+n}=[V_{2,m-1}]_{j,j+n}=0$ and $a^2+b^2=1$, $a^2-b^2=\rho$. We prove the following lemma.
\begin{lemma}\label{remain1}
Suppose the conditions of Theorem~\ref{shift} hold, we have
\begin{equation*}
\max\{|\mathcal A_{112}|,|\mathcal A_{114}|\}\le \frac C{nv^2}.
\end{equation*}
\end{lemma}
\begin{proof}
It is straightforward to check that
\begin{align*}
|\mathcal A_{112}|&\le \frac1{n^2}\E\|\V_{2,m}\J\RR\|_2^2,\\
|\mathcal A_{114}|&\le \frac1{n^2}\E^{\frac12}\|\V_{2,m}\J\RR\V_{2,m-1}\|_2^2
\E^{\frac12}\|\J\RR\|_2^2.
\end{align*}
Using well-known properties of Frobenius norm for matrices, $\|\mathbf A\mathbf B\|_2=\|\mathbf B\mathbf A\|_2$ and
$\|\mathbf A\mathbf B\|_2\le \|\mathbf A\|\|\mathbf B\|_2$, we get
\begin{align*}
|\mathcal A_{112}|&\le \frac1{n^2v^2}\E\|\V_{2,m}\|_2^2,\\
|\mathcal A_{114}|&\le
\frac1{n^2}\E^{\frac12}\|\V_{2,m}\J\RR\V_{1,m-1}\|_2^2\E^{\frac12}\|\J\RR\|_2^2.
\end{align*}
Furthermore, we note
\begin{align*}
\E\|\V_{2,m}\J\RR\V_{1,m-1}\|_2^2&=\E\|{\mathbf H^{(1)}}^{-1}\V_{1,m}\J\RR\V_{1,m-1}\|_2^2\\
&=\E\|{\HH^{(1)}}^{-1}(\I+\alpha\RR+\J(z)\RR)\V_{1,m-1}\|_2^2\\
&\le\E\Big(\|\V_{2,m-1}+{\HH^{(1)}}^{-1}(\alpha\I+\J(z))\RR\V_{1,m-1}\|_2^2\Big)\\
&\le 2\Big(\E\|\V_{2,m-1}\|_2^2+\frac{(|\alpha|+|z|)^2}{v^2}\E\|\V_{1,m-1}{\HH^{(1)}}^{-1}\|_2^2\Big)\\
&\le 2(1+\frac{(|\alpha|+|z|)^2}{v^2})\E\|\V_{2,m-1}\|_2^2.
\end{align*}
Applying Lemma \ref{norm4}, we conclude the proof.
\end{proof}
By Lemma \ref{remain1}, we may write
\begin{equation*}
\mathcal A_{1}=\mathcal A_{111}+\mathcal A_{113}+\varepsilon_n(\alpha,z).
\end{equation*}

\begin{lemma}\label{a111}
Under the assumptions of Theorem \ref{shift} we have
\begin{equation*}
|\mathcal A_{111}|\le \frac C{nv^2} ( 1 + v^{-2}).
\end{equation*}
\end{lemma}
\begin{proof}  A simple calculation shows that
\begin{equation*}
I:= \frac{2}{n^2}\sum_{j=1}^{n-1}\sum_{k=j+1}^n\E[\V_{2,m}\J\RR]_{jj}[\V_{2,m}
\J\RR]_{kk}=\frac{1}{n^2}\sum_{1\le j\ne k\le n}\E[\V_{2,m}\J\RR]_{jj}[\V_{2,m}\J\RR]_{kk}
\end{equation*}
By Lemma \ref{remain1}
\begin{equation*}
\frac1{n^2}\sum_{j=1}^n\E[\V_{2,m}\J\RR]_{jj}^2\le \frac C{nv^2}.
\end{equation*}
We may write
\begin{equation*}
I=\E\left(\frac1n\sum_{j=1}^n[\V_{2,m}\J\RR]_{jj}\right)^2+\varepsilon_n(\alpha,z).
\end{equation*}
Applying  Lemma \ref{var1}, we obtain
\begin{equation}\label{super}
\left |I-\left(\frac1n\sum_{j=1}^n\E[\V_{2,m} \J\RR]_{j,j}\right)^2 \right |\le \frac C{nv^2}(1+v^{-2}).
\end{equation}
Note that $\HH^{(q)}$, $q=1,\ldots,m$ have a symmetric joint distribution of entries,
i.e. $\HH^{(q)}$ has the same joint distribution of entries as $-\HH^{(q)}$, for any $q=1,\ldots,m$.
It follows immediately that
\begin{equation}\label{symmetry}
 \E\Tr\V_{2,m}\J\RR=0.
\end{equation}
To prove~\eqref{symmetry} we may replace the matrices $\HH^{(1)}$ and $\mathbf H^{2)}$ in the definition of
$\V_{2,m}\J\RR$ by $-\HH^{(1)}$ and $-\HH^{(2)}$. The resolvent matrix $\mathbf R$ still the same, since
\newline $\prod_{q=1}^m\HH^{(q)}=(-\HH^{(1)})(-\HH^{(2)})\prod_{q=3}^m\HH^{(q)}$ and we get
\begin{equation}\notag
 \E[\V_{2,m}\J\RR]_{jj}=-\E[\V_{2,m}\J\RR]_{jj}=0.
\end{equation}
The inequality~\eqref{super}) and equality~\eqref{symmetry} together imply the result of Lemma. Thus Lemma~\ref{a111} is proved.
\end{proof}
Finally, we prove that
\begin{align*}
\mathcal A_1=&-\frac {2}{n^2}\sum_{j=1}^{n-1}\sum_{k=j+1}^n\E[\V_{2,m}\J\RR\V_{1,m-1}]_{k,k+n}[\J\RR]_{j+n,j}\\
&-\frac{2}{n^2}\sum_{j=1}^{n-1}\sum_{k=j+1}^n\E[\V_{2,m}\J\RR\V_{1,m-1}]_{j,j+n}[\J\RR]_{k+n,k}+\varepsilon_n(\alpha,z).
\end{align*}
This equality we may rewrite as follows
\begin{align}\label{a1o}
\mathcal A_1&=-\frac {1}{n^2}\E\sum_{1\le j\ne k\le n}([\V_{2,m}\J\RR\V_{1,m-1}]_{k,k+n}
[\mathbf J\mathbf R]_{j+n,j}\nonumber\\
&\qquad\qquad\qquad+\E[\V_{2,m}\J\RR\V_{1,m-1}]_{j,j+n}[\J\RR]_{k+n,k})+\varepsilon_n(\alpha,z).
\end{align}
It is straightforward to check that
\begin{align}\label{a1b}
&\frac {1}{n^2}|\E\sum_{j=1}^{n}[\V_{2,m}\J\RR\V_{1,m-1}]_{j,j+n}[\J\RR]_{j+n,j}|\nonumber\\
&\qquad\qquad\qquad\qquad\qquad\le \frac C{n^{\frac32}v}\E^{\frac12}\|\V_{2,m}\J\RR\V_{1,m-1}\|_2^2\le \frac C{ nv^2}.
\end{align}
Relations (\ref{a1o}) and (\ref{a1b}) together imply
\begin{align*}
\mathcal A_1&=-\frac {1}{n^2}\E\sum_{j=1}^n\sum_{k=1}^n([\V_{2,m}\J\RR\V_{1,m-1}]_{k,k+n}
[\J\RR]_{j+n,j}\\
&\qquad\qquad\qquad+\E[\V_{2,m}\J\RR\V_{1,m-1}]_{j,j+n}[\J\RR]_{k+n,k})+\varepsilon_n(\alpha,z).
\end{align*}
By Lemma~\ref{var0}, Lemma~\ref{var1} and $\frac1n\sum_{j=1}^n\E[\J\RR]_{j,j+n}=
\frac1n\sum_{j=1}^n\E[\J\RR]_{j+n,j}=s_n(\alpha,z)$, we get
\begin{align}\label{A1}
\mathcal A_1&=-s_n(\alpha,z)\frac1n\sum_{j=1}^n\E([\V_{2,m}\J\RR\V_{1,m-1}]_{j,j+n}\nonumber \\
&\qquad\qquad\qquad\qquad+[\V_{2,m}\J\RR\V_{1,m-1}]_{j+n,j})+\varepsilon_n(\alpha,z).
\end{align}
Consider now the quantity $\mathcal A_2$. Similar to (\ref{A1}), we obtain
\begin{align*}
%\label{A2}
\mathcal A_2&=-s_n(\alpha,z)(\frac1n\sum_{j=1}^n\E([\V_{2,m}\J\RR\V_{1,m-1}]_{j+n,j}\\
&\qquad\qquad\qquad\qquad\qquad+[\V_{2,m}\J\RR\V_{1,m-1}]_{j,j+n})+\varepsilon_n(\alpha,z).
\end{align*}
Introduce the  notation, for $\nu=2,\ldots, m$
\begin{equation}\label{fq1}
 f_{q}=\frac1n\sum_{j=1}^n\E[\V_{q,m}\J\RR\V_{1,m-q+1}]_{j,j+n}
\end{equation}
and
\begin{equation*}
f_{m+1}=\frac1n\sum_{j=1}^n\E[\J\RR]_{j,j+n}=s_n(\alpha,z).
\end{equation*}
We rewrite the equality~\eqref{A1} using these notations
\begin{equation}\label{super1}
\mathcal A_1=-f_2s_n(\alpha,z)+\varepsilon_n(\alpha,z).
\end{equation}
We shall investigate the asymptotic of $f_{q}$, for $q=2,\ldots,m$. By definition of the matrices $\V_{q,m}$ and $\HH^{(q)}$, we have
\begin{equation*}
 f_{q}=\frac1{n\sqrt n}\sum_{k,j=1}^n\E Y^{(q)}_{jk}[\V_{q+1,m}\J\RR\V_{1,m-q+1}]_{k,j+n}.
\end{equation*}
We represent $f_q$ in the form
\begin{equation}\label{a1}
f_q=f_{q1}+f_{q2}+f_{q3},
\end{equation}
where
\begin{align*}
&f_{q1}=\frac1{n\sqrt n}\sum_{j=1}^{n-1}\sum_{k=j+1}^n\E Y^{(q)}_{jk}[\V_{q+1,m}\J\RR\V_{1,m-q+1}]_{k,j+n},\\
&f_{q2}=\frac1{n\sqrt n}\sum_{j=1}^{n}\E Y^{(q)}_{jj}[\V_{q+1,m}\J\RR\V_{1,m-q+1}]_{j,j+n},\\
&f_{q3}=\frac1{n\sqrt n}\sum_{j=2}^{n}\sum_{k=1}^{j-1}\E Y^{(q)}_{jk}[\V_{q+1,m}\J\RR\V_{1,m-q+1}]_{k,j+n}.
\end{align*}
Similarly to the previous steps we get
\begin{align*}
f_{q}&=\frac1n\sum_{k=1}^n\E[\V_{q+1,m}\J\RR\V_{1,m-q}]_{k,k+n}\\
&-\frac1n\sum_{k=1}^n\E[\mathbf V_{q+1,m}\J\RR\V_{1,m-q}]_{k,k+n}\frac1n\sum_{j=1}^n\E[\V_{m-q+2,m}\J
\RR\V_{1,m-q+1}]_{j+n,j+n}\\
&=f_{\nu+1}(1-\frac1n\sum_{j=1}^n\E[\V_{m-q+2,m}\J\RR\V_{1,m-q+1}]_{j+n,j+n})+\varepsilon_n(\alpha,z).
\end{align*}
Note that
\begin{align*}
\frac1n\sum_{j=1}^n\E[\V_{m-\nu+2,m}\J\RR\V_{1,m-\nu+1}]_{j+n,j+n}= \frac1n\sum_{j=1}^n\E[\V_{1,m}\J\RR]_{j+n,j+n}.
\end{align*}
Furthermore,
\begin{equation}\label{a10}
\frac1n\sum_{j=1}^n\E[\V_{1, m}\J\RR]_{j+n,j+n}=1+\alpha s_n(\alpha,z)+\overline z u_n(\alpha,z).
\end{equation}
Relations~\eqref{fq1}--\eqref{a10} together imply
\begin{equation*}
f_{q}=f_{q+1}(-\alpha s_n(\alpha,z)-\overline zu_n(\alpha,z))+\varepsilon_n(\alpha,z).
\end{equation*}
By induction we get
\begin{equation}\label{super100}
f_2=(-1)^{m-1}(\alpha s_n(\alpha,z)+\overline z u_n(\alpha,z))^{m-1}s_n(\alpha,z)+\varepsilon_n(\alpha,z).
\end{equation}
Relations (\ref{super1}) and (\ref{super100}) together imply
\begin{equation}\label{super101}
\mathcal A_1=(-1)^{m}(\alpha s_n(\alpha,z)+\overline z u_n(\alpha,z))^{m-1}s_n^2(\alpha,z)+\varepsilon_n(\alpha,z).
\end{equation}
Introduce now the  notations
\begin{equation*}
 h_{q}=\frac1n\sum_{j=1}^n\E[\V_{q,m}\J\RR\V_{1,m-q+1}]_{j+n,j},
\end{equation*}
for $q=2,\ldots, m$, and
\begin{equation*}
 h_{m+1}=\frac1n\sum_{j=1}^n\E[\J\RR]_{j+n,j}=s_n(\alpha,z).
\end{equation*}
Similar to~\eqref{super100} we get that
\begin{equation}\label{super102}
{h}_2=(-1)^{m-1}(\alpha s_n(\alpha,z)+z t_n(\alpha,z))^{m-1}s_n(\alpha,z)+\varepsilon_n(\alpha,z).
\end{equation}
and
\begin{equation}\label{b2}
\mathcal A_2=(-1)^{m}(\alpha s_n(\alpha,z)+z t_n(\alpha,z))^{m-1}s_n^2(\alpha,z)+\varepsilon_n(\alpha,z).
\end{equation}
Consider now the function $t_n(\alpha,z)$ which we may represent  as follows
\begin{equation*}
\alpha t_n(\alpha,z)=\frac1n\sum_{j=1}^n\E[\mathbf V(z)\mathbf R]_{j+n,j}.
\end{equation*}
By definition of the matrix $\mathbf H^{(1)}$, we may write
\begin{equation}\label{tn1}
\alpha t_n(\alpha,z)=\frac1n\sum_{j,k=1}^n\E Y_{jk}^{(m)}[\V_{2,m}\J\RR]_{j+n,k}-\overline z\ s_n(\alpha,z).
\end{equation}
The first term in the r.h.s. of (\ref{tn1}) we represent in the form
\begin{equation*}
\mathcal B_1:=\frac1n\sum_{j,k=1}^n\E Y_{jk}^{(m)}[\V_{2,m}\J\RR]_{j+n,k}=\mathcal B_{11}+\mathcal B_{12}+\mathcal B_{13},
\end{equation*}
where
\begin{align*}
\mathcal B_{11}&=\frac1n\sum_{j=1}^{n-1}\sum_{k=j+1}^n\E Y_{jk}^{(m)}[\V_{2,m}\J\RR]_{j+n,k},\\
\mathcal B_{11}&=\frac1n\sum_{j=1}^{n}\E Y_{jk}^{(m)}[\V_{2,m}\J\RR]_{j+n,k},\\
\mathcal B_{13}&=\frac1n\sum_{j=2}^{n}\sum_{k=1}^{j-1}\E Y_{jk}^{(m)}[\V_{2,m}\J\RR]_{j+n,k},
\end{align*}
Previous relations together imply
\begin{align*}
\alpha t_n(\alpha,z)&=-\frac1n\sum_{j=1}^n\E[\V_{2,m}\J\RR\V_{1,m-1}]_{j+n,j}\frac1n\sum_{k=1}^n\E[\RR]_{k+n,k}\\
&-\overline z s_n(\alpha,z)+\varepsilon_n(\alpha,z)\\
&={h}_2\ t_n(\alpha,z)-
\overline z\ s_n(\alpha,z)+\varepsilon_n(\alpha,z).
\end{align*}
Applying the equality (\ref{super102}), we obtain
\begin{align}\label{super105}
\alpha t_n(\alpha,z)=(-1)^{m}(\alpha s_n(\alpha,z)+z t_n(\alpha,z))^{m-1}&s_n(\alpha,z)t_n(\alpha,z)\nonumber\\
&-\overline z\ s_n(\alpha,z)+\varepsilon_n(\alpha,z).
\end{align}
Analogously we obtain
\begin{align}\label{super106}
\alpha u_n(\alpha,z)=(-1)^{m}(\alpha s_n(\alpha,z)+\overline z u_n(\alpha,z))^{m-1}&s_n(\alpha,z)u_n(\alpha,z)\nonumber\\
&-z\ s_n(\alpha,z)+\varepsilon_n(\alpha,z).
\end{align}
Since $|\alpha|\ge v$, we may rewrite these equation as follows
\begin{align}
t_n(\alpha,z)&=(-1)^{m}(\alpha s_n(\alpha,z)+z t_n(\alpha,z))^{m-1}\alpha^{-1}s_n(\alpha,z)t_n(\alpha,z)\nonumber\\
&\qquad\qquad\qquad\qquad\qquad-\overline z\ s_n(\alpha,z)\alpha^{-1}+\varepsilon_n(\alpha,z)\notag\\
u_n(\alpha,z)&=(-1)^{m}(\alpha s_n(\alpha,z)+\overline z u_n(\alpha,z))^{m-1}\alpha^{-1}s_n(\alpha,z)u_n(\alpha,z)\nonumber\\
&\qquad\qquad\qquad\qquad\qquad-z\ s_n(\alpha,z)\alpha^{-1}+\varepsilon_n(\alpha,z).\notag
\end{align}
The rest of the proof is the same as in the proof of Theorem 3.1, \cite{GotTikh2011}, p. 11-13. For the readers convenience we repeat it here. We note that, for some numerical constant $C>0$,
\begin{align}\label{in-1}
|\alpha s_n(\alpha,z)|\le 1+\Big|\frac1{2n}\E\Tr\mathbf R\mathbf V\Big|&\le 1+v^{-1}\frac Cn(\E^{\frac12}\|\mathbf W\|_2+n|z|)\nonumber\\
&\le C(1+\frac {|z|}v),
\end{align}
and
\begin{equation}\label{in-2}
\max\{|\overline zt_n(\alpha,z)|,| z u_n(\alpha,z)|\}\le \frac{|z|}v.
\end{equation}
Introduce notation
\begin{align*}
P:=P(\alpha,z)&=\alpha s_n(\alpha,z)+\overline z u_n(\alpha,z)\\
Q:=Q(\alpha,z)&=\alpha s_n(\alpha,z)+z t_n(\alpha,z).
\end{align*}
Multiplying~\eqref{super105} by $z$ and~\eqref{super106} by $\overline z$ and subtracting the second  one from the first equation, we obtain
\begin{align}\label{in-3}
zt_n(\alpha,z)-\overline zu_n(\alpha,z)&=(zt_n(\alpha,z)-\overline zu_n(\alpha,z)) \nonumber\\
&\times s_n(\alpha,z)zt_n(\alpha,z)\alpha^{-1}(P^{m-2}+QP^{m-3}+\cdots+Q^{m-2})\nonumber\\
&+Q^{m-1}s_n(\alpha,z)\alpha^{-1}( zt_n(\alpha,z)-\overline zu_n(\alpha,z))+\varepsilon(\alpha,z).
\end{align}
Using inequalities (\ref{in-1}), (\ref{in-2}) and $|s_n(\alpha,z)|\le v^{-1}$, we get
\begin{align}\label{in-4}
|s_n(\alpha,z)zt_n(\alpha,z)\alpha^{-1}
(P^{m-2}+QP^{m-3}+\cdots+Q^{m-2})|&\le \frac{C^{m-1}m(1+\frac{|z|}v)^{m-2}}{v^3},\nonumber\\
|Q^{m-1}s_n(\alpha,z)\alpha^{-1}|\le \frac{C^{m-1}(1+\frac{|z|}v)^{m-2}}{v^3}.
\end{align}
From relations (\ref{in-3}) and (\ref{in-4}) we may conclude that there exists $V_0=V_0(m,z)$ depending on $m$ and $z$ such that for all $v\ge V_0$ 
\begin{equation}\label{superchto}
z t_n(\alpha,z)=\overline z u_n(\alpha,z)+\varepsilon_n(\alpha,z).
\end{equation}
The last relation implies that
\begin{equation}\label{b3}
\mathcal A_1=\mathcal A_2+\varepsilon_n(\alpha,z).
\end{equation}
Relations (\ref{b1}), \ref{super101}), (\ref{b2}), (\ref{superchto}, and (\ref{b3}) together imply
\begin{align}
1+\alpha s_n(\alpha,z) =(-1)^m(\alpha s_n(\alpha,z)+z \ t_n(\alpha,z))^{m-1}&s_n^2(\alpha,z) \nonumber\\
&-z \ t_n(\alpha,z)+\varepsilon_n(\alpha,z). \label{supernu}
\end{align}
Introduce the notations
\begin{equation*}
g_n:=s_n(\alpha,z),\quad w_n:=\alpha+\frac{z\ t_n(\alpha,z)}{g_n}.
\end{equation*}
Using these notations  we may rewrite the equations (\ref{supernu}) and (\ref{superchto}) as follows
\begin{align}\label{system}
&1+w_ng_n=(-1)^{m}g_n^{m+1}w_n^{m-1}+\varepsilon_n(\alpha,z)\nonumber\\
&(w_n-\alpha)+(w_n-\alpha)^2g_n-g_n|z|^2=\varepsilon_n(\alpha,z).
\end{align}
Let $n,n'\to\infty$. Consider the  difference $g_n-g_{n'}$. From the first inequality it follows that
\begin{equation*}
|g_n-g_{n'}|\le \frac{|\varepsilon_{n,n'}(\alpha,z)|+|w_n-w_{n'}||g_n+(-1)^{m+1}g_{n'}^{m+1}
(w_n^{m-2}+\cdots+w_{n'}^{m-2})|}{|w_n+(-1)^{m+1}y_{n'}^{m+1}(w_n+(-1)^{m+1}w_{n}^{m-1}(g_n^{m}+
\cdots+g_{n'}^{m})|}
\end{equation*}
Note that $\max\{|g_n|,\ |g_{n'}|\}\le \frac1v$ and $\max\{|w_n|,\ |w_{n'}|\}\le C+v$ for some positive constant
$C=C(m)$ depending of $m$. We may choose a sufficiently large $V_0'$ such that for any $v\ge V_0'$ we obtain
\begin{equation}\label{100}
|g_n-g_{n'}|\le \frac{|\varepsilon_{n,n'}(\alpha,z)|}{v}+\frac Cv|w_n-w_{n'}|.
\end{equation}
Furthermore, the second equation {in (\ref{system})} implies that
\begin{align*}
&(w_n-w_{n'})(1+g_n(w_n+w_{n'}-2\alpha))\\
&\qquad\qquad\qquad\qquad=(g_n-g_{n'})((w_n-\alpha)^2-|z|^2)+\varepsilon_{n,n'}(\alpha,z).
\end{align*}
It is straightforward to check that $\max\{|w_n-\alpha|, |w_{n'}-\alpha|\}\le (1+|\varepsilon_n(\alpha,z)|)|z|$.
This implies that there exists $V_1$ such that for any $v\ge V_1$
\begin{equation}\label{101}
|w_n-w_{n'}|\le |\varepsilon_{n,n'}(\alpha,z)|+4|z|^2|g_n-g_n'|.
\end{equation}
Inequalities (\ref{100}) and (\ref{101}) together imply that there exists a constant $V_0=\max\{V_0',V_1\}$ such that
for any $v\ge V_0$
\begin{equation*}
|g_n-g_n'|\le |\varepsilon_{n,n'}(\alpha,z)|,
\end{equation*}
where $\varepsilon_{n,n'}(\alpha,z)\to 0$ as $n,{n}'\to \infty$ uniformly with respect to $v\ge V_0$ and
$|u|\le C$ ($\alpha=u+iv$). 

Since $g_n,g_{n'}$ are  locally bounded {analytic functions in the upper half-plane}
we may conclude by Montel's Theorem (see, for instance, \cite{Conw1978}, p. 153, Theorem 2.9) that there exists an analytic  function $g_0$   in the upper half-plane such that $\lim g_n=g_0$. Since $g_n$ are Nevanlinna functions, (that is analytic functions mapping the upper half-plane into itself)  $g_0$ will be a Nevanlinna function too and there exists some distribution function $G(a,z)$ such that 
$$
g_0=\int_{-\infty}^{\infty}\frac1{a-\alpha}dG(a,z)
$$ 
and
\begin{equation*}\label{pip}
\Delta_n(z):=\sup_a|G_n(a,z)-G(a,z)|\to 0\quad\text{as}\quad n\to\infty.
\end{equation*}

The function $g_0$ satisfies the equations (\ref{system00}). Thus Proposition \ref{respect}  is proved.
\end{proof}

The Lemma \ref{universality} and Proposition \ref{respect} together conclude the proof of Theorem \ref{shift}.
Thus Theorem \ref{shift} is proved.

\section{The minimal singular value of matrix $\mathbf V(z)$}\label{singular}
We shall use  the following theorem which was proved in~\cite{GotNauTikh2013sing}.
%%%%%%%%%%%%%%%%%%%%%%%%%%%%%%%%%%%%%%%%%%%%%%%%%%%%%%%%%%%%%%%%%%%%%%%%%%%%%%%%%%%%%%%%%%%%%%%%%%%%%%%%%%%%%%%%%%%%%%%%%%%%%%%%%%
%%%%%%%%%%%%%%%%%%%%%%%%%%%%%%%%%%%%%%%%%%%%%%%%%%%%%%%%%%%%%%%%%%%%%%%%%%%%%%%%%%%%
\begin{theorem} \label{thm10}
Assume that $X_{jk}$, $1\le j,k\le n$, satisfy the conditions \Cond and \UI.  Let $\X=\{X_{jk}\}$ denote a $n\times n$ random matrix with the entries $X_{jk}$ and let  $\mathbf M_n$ denote a non-random matrix with $\|\mathbf M_n\|\le Kn^Q=:K_n$ for some $K>0$  and $Q\ge0$. Then there exist constants $ C, A,B >0$ depending on  $K,Q$ and $\rho$ such that
\begin{align}
\Pb (s_n\le n^{-B } )\le Cn^{-A},
\end{align}
\end{theorem}

\begin{lemma}\label{spec}
 Under the conditions of Theorem \ref{th:main} there exists a constant $C$ such that for any   $k\le n(1-
C\Delta_n^{\frac1{m+1}}(z))$,
\begin{equation}\notag
 \Pb\{s_k\le \Delta_n(z)\}\le C\Delta_n^{\frac1{m+1}}(z).
\end{equation}
\end{lemma}
\begin{proof}
We may write, for any $k=1,\ldots,n$,
{\begin{equation}\notag
 \Pb\{s_k\le \Delta_n(z)\}\le \Pb\{{\overline{\mathcal G}_n}(s_k,z)\le {\overline{\mathcal G}_n (\Delta_n(z),z)}\}\le
 \Pb\{\frac{n-k}n\le{ \overline{\mathcal G}_n(\Delta_n(z),z)}\}.
\end{equation}}
Applying  Chebyshev's  inequality, we obtain
\begin{equation*}
 \Pb\{s_k\le \Delta_n(z)\}\le \frac{n\E\overline{\mathcal G}_n(\Delta_n(z),z)}{n-k}\le
\frac{n(\overline G(\Delta_n(z),z)+2\Delta_n(z)}{n-k}.
\end{equation*}
It is straightforward to check that from the system of equations  (\ref{system00}) it follows
\begin{equation*}
 \overline G(\Delta_n(z),z)\le C\Delta_n^{\frac2{m+1}}(z).
\end{equation*}
The last inequality concludes the proof of Lemma \ref{spec}.
\end{proof}

\begin{lemma}\label{l:small sing val}
Let $n_1:=[n-n\delta_n]+1$ and $n_2:=[n-n^{\gamma}]$ for any sequence $\delta_n\to0$, and some  $0<\gamma<1$. Under the conditions of Theorem~\ref{th:main} we have
\begin{align}
&\lim_{n\to\infty}\frac1n\sum_{n_1\le j\le n_2}\ln s_j(\mathbf X^{(q)})=0,\quad\text{for}
\quad q=1,\ldots,m-1,\notag\\
&\lim_{n\to\infty}\frac1n\sum_{{n_1}\le j\le {n_2}}\ln s_j(\mathbf X^{(m)}+\M_n)=0\notag,
\end{align}
where $||\M_n|| \le n^{Q}$ for some $Q > 0$.
\end{lemma}
\begin{proof}
The claim follows from the bound
\begin{equation}\label{prod3}
s_j(\mathbf X^{(\nu)}+\mathbf M_n)\ge c \frac{n-j}{n},\quad 1\le j\le n-n^{\gamma}.
\end{equation}
To prove this we need  the following simple Lemma.
\begin{lemma}\label{vot}
 Let $\lim_{n\to\infty}\delta_n=0$ and let  $q_j$, for $n_1\le j\le
 n_2$ with $0<\gamma<1$  denote numbers
 satisfying  the inequalities
\begin{equation}\notag
 n^Q\ge q_j\ge c \frac{n-j}n
\end{equation}
for some constant $Q>0$.
Then
\begin{equation}\notag
 \lim_{n\to\infty}\frac1n\sum_{n_1\le j\le n_2}\ln q_j=0.
\end{equation}
\end{lemma}

\begin{proof}
Note that
\begin{equation}\notag
0\le\frac1n\sum_{n_1\le j\le n_2: \ q_j\ge 1}\ln q_j\le Qn^{-(1-\gamma)}\ln n\to0, \text{ as }n\to\infty.
\end{equation}
Without loss of generality we may assume that $0<q_j\le 1$. By the conditions of {Lemma \ref{vot}}, we have
\begin{equation*}
 0\ge\frac1n\sum_{n_1\le j\le n_2}\ln q_j\ge\frac1n\sum_{n_1\le j\le n_2}
\ln \{\frac{n-j}n\}=A.
\end{equation*}
After summation  and using Stirling's formula, we get
\begin{align}\label{009}
|A|\le \frac1n\ln\Big\{\frac{n_1!}{n_2!n^{n_2-n_1}}\Big\}\le
\delta_n|\ln\delta_n|+(1-\gamma)n^{\gamma-1}\ln n\to 0 \text{ as } n \to \infty.
\end{align}
This  proves  Lemma \ref{vot}.
\end{proof}
We continue the proof of Lemma~\ref{l:small sing val}. It remains to prove the inequality (\ref{prod3}).
Similar  result for matrices with independent entries  was proved by Tao and Vu in~\cite{TaoVu2010}
(see inequality (8.4) in~\cite{TaoVu2010}). It represents the  crucial result in their proof of the
circular law assuming the second moment only. For completeness  we give here a simple modification of their proof for the case of random matrices with correlated entries. We start from the following
\begin{statement}\label{Tao}
Let $1\le d\le n-n^{\gamma}$ with $\frac{8}{15}<\gamma<1$.
and $0<c<1$, and $\mathbb H$ be a (deterministic) $d$-dimensional subspace of $\mathbb C^n$. Let $X_{j}$
be independent  random variables with $\E X_{j}=0$ and $\E|X_{j}|^2=1$, squares of
which are uniformly integrable , i.e.
\begin{equation}\label{unif0+}
\max_{j}\E|X_{j}|^2\mathbb I{\{|X_{j}|>M\}}\to0\quad\text{as}\quad M\to \infty.
\end{equation}
Let $\mathbf x^T=(X_1,\ldots,X_n)+(m_1,\ldots,m_N)$ where $\mathbf m^T=(m_1+\ldots,m_n)$ is non-random vector. Then
\begin{equation}\label{dist100}
\Pb\{\text{\rm dist}(\mathbf x+\mathbf m,\mathbb H)\le c\sqrt{n-d}\}=O(\exp\{-n^{\frac{\gamma}8}\}),
\end{equation}
where $\text{\rm dist}(X,\mathbb H)$ denotes the Euclidean distance between a vector $X$ and a subspace
$\mathbb H$ in $\mathbb C^n$.
\end{statement}
\begin{proof} It was proved by Tao and Vu in~\cite{TaoVu2010} (see Proposition 5.1).
Here we sketch their proof. As shown in~\cite{TaoVu2010} we may reduce the problem to the case that $\E X=0$.
For this it is enough to consider vectors $\mathbf x'$ and $\mathbf v$  such that $\mathbf x=\mathbf x'+\mathbf v$ and $\E \mathbf x'=0$.
Instead of the subspace $\mathbb H$ we may consider subspace $\mathbb H'=\text{\rm span}(\mathbb H,\mathbf v)$ and note that
\begin{equation}
 \text{\rm dist}(\mathbf x,\mathbb H)\ge \text{\rm dist}(\mathbf x',\mathbb H').
\end{equation}
The claim follows now from a  corresponding result for random vectors with mean zero. In what follows we assume that $\E \mathbf x=0$. We reduce the problem to vectors with  bounded coordinates.
Let $\xi_j=\mathbb I\{|X_j|\ge n^{\frac{1-\gamma}2}\}$, where $X_j$ denotes the $j$-th coordinate of a vector $\mathbf x$. Note that $p_n:=\E\xi_j\le n^{-(1-\gamma)}$. Applying Chebyshev's  inequality, we get, for any $h>0$
\begin{equation}\notag
\Pb\{\sum_{j=1}^n\xi_j\ge 2{n^{\gamma}}\}\le \exp\{-hn^{\gamma}\}\exp\{np_n(\text{\rm e}^h-1-h)\}.
\end{equation}
Choosing $h=\frac14$, we obtain
\begin{equation}\label{prod7}
\Pb\{\sum_{j=1}^n\xi_j\ge 2{n^{\gamma}}\}\le\exp\{-\frac{n^{\gamma}}8\}.
\end{equation}
Let $J\subset\{1,\ldots,n\}$ and $E_J:=\{\prod_{j\in J}(1-\xi_j)\prod_{j\notin J}\xi_j=1\}$.
Inequality (\ref{prod7}) implies
\begin{equation}\notag
 \Pb\{\bigcup_{J:|J|\ge n-2n^{\gamma}}E_J\}\ge 1-\exp\{-\frac{n^{\gamma}}8\}.
\end{equation}
Let $J$ with $|J|\ge n-2n^{\gamma}$ be fixed. Without loss of generality we may assume that $J=1,\ldots,n'$
with some $n-2n^{\gamma}\le n'\le n$. It is now sufficient to prove that
\begin{equation}\label{dist1}
 \Pr\{\text{\rm dist}(\mathbf x,\mathbb H)\le c\sqrt{n-d}| E_J\}=O(\exp\{-\frac {n^{\gamma}}8\}).
\end{equation}
Let $\pi$ denote the orthogonal projection $\pi:\mathbb C^n\rightarrow\mathbb C^{n'}$. We note that
\begin{equation}\label{dist9}
\text{\rm dist}(\mathbf x,\mathbb H)\ge\text{\rm dist}(\pi(\mathbf x),\pi(\mathbb H)).
\end{equation}

Let $\widetilde X$ be a random variable $X$ conditioned on
 the event $|X|\le n^{1-\gamma}$ and let
$\widetilde {\mathbf x}=(\widetilde x_1,\ldots, \widetilde x_n)$. The relation (\ref{dist1}) will follow now from
\begin{equation}\notag
 \Pb\{\text{\rm dist}(\widetilde{\mathbf x'},\mathbb H')\le c\sqrt{n-d}\,\big| |x_j|\le n^{1-\gamma}, j\notin J\}
=O(\exp\{-\frac {n^{\gamma}}8\}),
\end{equation}
where $\mathbb H'=\pi(\mathbb H)$ and $\widetilde {\mathbf x'}=\pi(\widetilde {\mathbf x})$.
We may represent the vector $\widetilde {\mathbf x}$as  $\widetilde {\mathbf x}=\widetilde {\mathbf x'}+\mathbf v$, where
$\mathbf v=\E\widetilde {\mathbf x}$ and $\E\widetilde {\mathbf x'}=0$.
We reduce the claim to the bound
\begin{equation}\label{dist5}
\Pb\{\text{\rm dist}(\widetilde {\mathbf x'},\mathbb H'')\le c\sqrt{n-d}\,\big| |x_j|\le n^{1-\gamma}, j\notin J\}=
O(\exp\{-\frac {n^{\gamma}}8\}),
\end{equation}
where $\mathbb H''=\text{\rm span}(v,\mathbb H')$.
In what follows we shall omit the  symbol $'$ in the notations. To prove (\ref{dist5}) we shall apply the following result of Maurey. Let $\mathbb X$ denote a normed space and $f$ denote a convex function on $\mathbb X$.
Define the functional $Q$ as follows
\begin{equation}\notag
 Qf(x):=\inf_{y\in \mathbb X}[f(y)+\frac{\|x-y\|^2}4].
\end{equation}
\begin{definition}
 We say that a measure $\mu$ satisfies the convex property $(\tau)$ if for any convex function $f$ on $\mathbb X$
\begin{equation}\notag
 \int_{\mathbb X}\exp\{Qf\}d\mu\int_{\mathbb X}\exp\{-f\}d\mu\le 1.
\end{equation}

\end{definition}

We reformulate the following result of Maurey (see \cite{Maur1991}, Theorem 3).
{Following Maurey we shall say that $\nu$ has diameter $\le 1$ as a short way to express that $\nu$ is supported by a set of
diameter $\le1$.}
\begin{theorem}\label{maurey}Let $(\mathbb X_i)$ be a family of normed spaces; for each $i$,
let $\nu_i$ be a probability measure with diameter $\le 1$ on $\mathbb X_i$. If $\nu$ is the product of a family $(\nu_i)$, then $\nu$ satisfies the convex property $(\tau)$.
\end{theorem}
As corollary of Theorem \ref{maurey} we get
\begin{corollary}\label{talagrand}
  Let $\nu_i$ be a probability measure with diameter $\le 1$ on $\mathbb X$, $i=1,\ldots,n$. Let $g$ denote a
convex $1$-Lipshitz function on $\mathbb X^n$. Let $M(g)$ denote a median of $g$.
If $\nu$ is the product of the family $(\nu_i)$, then
\begin{equation}\notag
 \nu\{|g-M(g)|\ge h\}\le 4\exp\{-\frac {h^2}4\}.
\end{equation}
\end{corollary}
Applying Corollary \ref{talagrand} to $\nu_i$, being  the distribution of $\widetilde x_i$, we get
\begin{equation}\label{dist10}
\Pb\left\{|\text{\rm dist}(\widetilde {\mathbf x},\mathbb H)-M(\text{\rm dist}(\widetilde {\mathbf x},\mathbb H))|
\ge rn^{\frac{1-\gamma}2}\right\}\le 4\exp\{-r^2/16\}.
\end{equation}
The last inequality implies that there exists a constant $C>0$ such that
\begin{equation}\label{dist13}
 |\E\text{\rm dist}(\widetilde {\mathbf x},\mathbb H)-M(\text{\rm dist}(\widetilde {\mathbf x},\mathbb H))|\le Cn^{\frac{1-\gamma}2},
\end{equation}
and
\begin{equation}\label{dist14}
\E\text{\rm dist}(\widetilde {\mathbf x},\mathbb H)\ge \sqrt{\E(\text{\rm dist}(\widetilde {\mathbf x},\mathbb H))^2}-
Cn^{\frac{1-\gamma}2}.
\end{equation}
By Lemma 5.3 in~\cite{TaoVu2010}
\begin{equation}\label{dist15}
 \E(\text{\rm dist}(\widetilde {\mathbf x},\mathbb H))^2=(1-o(1){)}(n-d).
\end{equation}
Since $n-d\ge n^{\gamma}$ the inequalities (\ref{dist13}), (\ref{dist14}) and (\ref{dist15}) together imply
(\ref{dist100}).

Now we prove (\ref{prod3}). We repeat the proof of Tao and Vu~\cite{TaoVu2010}, inequality (8.4).
Fix $j$. Let $\mathbf A_n=\mathbf X^{(\nu)}-z\mathbf M_n$ and let $\mathbf A_n'$
denote a matrix formed  by the  first {$n'=n-k$} rows of $\sqrt{n}\mathbf A_n$ with $k=j/2$. Let {$\sigma_l$} ($\sigma_l'$), $1\le l\le n-k$,  be the singular values of {$\mathbf A_n$} ($\mathbf A_n')$ (in decreasing order). By the interlacing property and re-normalizing we get
\begin{equation}\notag
 \sigma_{n-j}\ge \frac1{\sqrt n}\sigma_{n-j}'.
\end{equation}
By Lemma A.4 in~\cite{TaoVu2010}
\begin{equation}\label{eq: T}
 T:={\sigma'_1}^{-2}+\cdots+{\sigma'_{n-k}}^{-2}={\text{\rm dist}}_1^{-2}+\cdots+{\text{\rm dist}}^{-2}_{n-k},
\end{equation}
{with
$$
{\text{\rm dist}}_{j}={\text{\rm dist}}(\mathbf x_j,\mathbb H_j),
$$
where $\mathbf x_j$ is the $j$-th row of matrix $\mathbf A_n'$ and $\mathbb H_j$ denotes hyperplane generated by the $n'-1$ rows $X_1,\ldots,X_{j-1},X_{j+1},\ldots,X_{n'}$. }
Let $\pi_j$ denote the projector onto $\mathbb R_j^{n-1}$ in $\mathbb R^{n}$
defined by $\pi_j(\mathbf x)=(X_1,\ldots,X_{j-1},0.X_{j+1},\ldots,X_n)$.
Then we have
\begin{equation}\notag
 {\text{\rm dist}}(\mathbf x_j,\mathbb H_j)\ge {\text{\rm dist}}(\pi_j(\mathbf x),\pi_j(\mathbb H_j)).
\end{equation}
Note that vector $\pi_j(\mathbf x)$ and subspace $\pi_j(\mathbb H_j)$ are independent and vector $\pi_j(\mathbf x)$ has independent coordinates. From~\eqref{eq: T}
\begin{equation}\notag
 T\ge (j-k){\sigma'}_{n-j}^{-2}=\frac j2{\sigma'}_{n-j}^{-2} \geq \frac{j}{2n}\sigma_{n-j}^{-2}.
\end{equation}
Applying Proposition \ref{Tao}, we get that with probability $1-\exp\{-n^{\gamma}\}$
\begin{equation}\notag
 T\le \frac{n}j.
\end{equation}
Combining the last inequalities, we get (\ref{prod3}). Thus Proposition \ref{Tao} is proved.
\end{proof}
This finishes the proof of Lemma.
\end{proof}

\begin{lemma}\label{l: log uniform integr}
Assume the assumptions of Theorem~\ref{th:main} hold, then $\ln( \cdot )$ is uniformly integrable in probability with respect to $\{\nu_n\}_{n \geq 1}.$
\end{lemma}

\begin{proof}[Proof of Lemma~\ref{l: log uniform integr}]
It is enough to check that
\begin{equation} \label{eq:log u.i.}
\lim_{t \rightarrow \infty} \varlimsup_{n \rightarrow \infty} \mathbb P \left (\int_0^\infty |\ln x| \nu_n(dx) > t \right ) = 0
\end{equation}

Let $k_0 = [n(1-C\Delta_n^{\frac1{m+1}}(z))]$. We introduce the event
\begin{align*}
\Omega_0:=\Omega_{0,n}:=\{\omega \in \Omega: &s_n(\X^{(q)}) \geq n^{-b}, q = 1, ... , m-1, \\
&s_n(\X^{(m)}+\M_n) \geq n^{-b}, s_{k_0} \geq \Delta_n(z)\}.
\end{align*}
for some $b > 0$ which will be chosen later and $\M_n = -z (\prod_{i=1}^{m-1} \X^{(q)})^{-1}$.
Note that the matrices $\mathbf X^{(m)}$ and $\mathbf M_n$ are independent and it follows from Theorem~\ref{thm10} that
$\|\mathbf M_n\|_2\le n^{ Q}$ for some ${Q}>0$ with probability close to one. From Theorem~\ref{thm10} and Lemma~\ref{l:small sing val} we conclude that $\varlimsup_{n \rightarrow \infty}\Pb(\Omega_0^{c}) = 0$.
It follows that it is enough to prove that
$$
\lim_{t \rightarrow \infty} \varlimsup_{n \rightarrow \infty} \mathbb P \left (\int_0^\infty |\ln x| \nu_n(dx) > t, \Omega_0 \right ) = 0
$$
We may split the integral $\int_0^\infty |\ln x| \nu_n(dx)$ into three terms
\begin{align*}
&T_1: = -\int_0^{\Delta_n} \ln x \nu_n(dx, z),\\
&T_2: = \int_{\Delta_n}^{\Delta_n^{-1}} |\ln x| \nu_n(dx, z),\\
&T_3: = \int_{\Delta_n^{-1}}^\infty \ln x \nu_n(dx, z).
\end{align*}
Denote by $n':=k_0+1$ and $n'':=[n - n^{1-\gamma}]$. We consider the term $T_1$ which we may rewrite as
$$
T_1 = -\frac{1}{n} \sum_{i = n' + 1}^n \ln s_i.
$$
We shall use the following well-known fact.
Let $\mathbf A$ and $\mathbf B$ be $n\times n$ matrices and let $s_1(\mathbf A)\ge\cdots \ge s_n(\mathbf A)$ resp.
($s_1(\mathbf B)\ge\cdots\ge s_n(\mathbf B)$ and
$s_1(\mathbf A\mathbf B)\ge\cdots \ge s_n(\mathbf A\mathbf B)$) denote
 the singular value of a matrix
$\mathbf A$ (and the matrices  $\mathbf B$ and  $\mathbf A\mathbf B$ respectively).
Then  we have
\begin{equation}\label{prod1}
 \prod_{j=k}^ns_j(\mathbf A\mathbf B)\ge \prod_{j=k}^ns_j(\mathbf A)s_j(\mathbf B),
\end{equation}
and
\begin{equation}\notag
\prod_{j=1}^ns_j(\mathbf A\mathbf B)= \prod_{j=1}^ns_j(\mathbf A)s_j(\mathbf B),
\end{equation}
for any $1\le k\le n$ (see, for instance \cite{HornJohn1990}, p.171, Theorem 3.3.4).
From~\eqref{prod1} it follows that
\begin{align*}
T_1 &\le -\frac{1}{n}\sum_{q = 1}^{m-1} \sum_{i = n'+1}^n \ln s_i(\X^{(q)}) - \frac{1}{n}\sum_{i = n'+1}^n \ln s_i(\X^{(m)}+\M_n) = \\
& -\frac{1}{n}\sum_{q = 1}^{m-1} \sum_{i = n'+1}^{n''} \ln s_i(\X^{(q)}) - \frac{1}{n}\sum_{i = n'+1}^{n''} \ln s_i(\X^{(m)}+\M_n)
 \\
& -\frac{1}{n}\sum_{q = 1}^{m-1} \sum_{i = n''+1}^{n} \ln s_i(\X^{(q)}) - \frac{1}{n}\sum_{i = n''+1}^{n} \ln s_i(\X^{(m)}+\M_n)
\end{align*}
From Lemma~\ref{l:small sing val}, inequality~\eqref{009} and definition of $\Omega_0$ it follows that
$$
T_1 \le C n^{\gamma-1} \ln n + \Delta_n |\ln \delta_n| \rightarrow 0 \text{ as } n \rightarrow \infty
$$
For the term $T_3$ we may write the bound
$$
T_3 \le \Delta_n |\ln \Delta_n| \int_0^\infty x^2 \nu_n(dx, z) \rightarrow 0 \text{ as } n \rightarrow \infty,
$$
where we have used the fact that $x^{-2} \ln x$ is a decreasing function for $x \geq \sqrt e$.
It remains to estimate $T_2$. Integrating by parts and using~\eqref{distance} we write
$$
\E T_2 \le C \Delta_n |\ln \Delta_n| + \int_{\Delta_n}^{\Delta_n^{-1}} |\ln x| dG(x, z) < \infty
$$
Using Markov's inequality we finish the proof of Lemma.
\end{proof}

\section{Appendix}

\begin{lemma}\label{mean}
 Under the conditions of Theorem \ref{th:main} we have, for any {$j,k=1,\ldots,n$, }
 and for any $1\le \alpha\le \beta\le m$,
\begin{equation}\notag
\E[\mathbf V_{\alpha,\beta}]_{jk}=0
\end{equation}
\end{lemma}
\begin{proof}For $\alpha=\beta$ the claim is easy. Let $\alpha<\beta$.
 Direct calculations show that
\begin{equation}\notag
 \E[\mathbf V_{\alpha,\beta}]_{jk}=\frac1{n^\frac{\beta-\alpha}2}\sum_{j_1=1}^{p_{\alpha}}
\sum_{j_2=1}^{p_{\alpha+1}}\dots \sum_{j_{\beta-\alpha}=1}^{p_{\beta-1}}
\E X^{(\alpha)}_{j,j_1}X^{(\alpha+1)}_{j_1,j_2}\cdots X^{(\beta)}_{j_{\beta-\alpha},k}=0
\end{equation}
Thus the Lemma is proved.
\end{proof}
In all Lemmas below we shall assume that
\begin{equation}\label{as1}
 \E X_{jk}^{(\nu)}=0,\quad\E|X_{jk}^{(\nu)}|^2=1, \quad |X_{jk}^{(\nu)}|\le c\tau_n\sqrt n\quad\text{a. s.}
\end{equation}

\begin{lemma}\label{norm2}
Under the conditions of Theorem \ref{th:main} assuming (\ref{as1}), we have,  for any $1\le \alpha\le \beta\le m$,
\begin{equation*}
\E\|\mathbf V_{\alpha,\beta}\|_2^2\le Cn
\end{equation*}
\end{lemma}
\begin{proof}We shall consider the case $\alpha<\beta$ only. The other cases are obvious. Direct calculation shows that
\begin{equation*}
\E\|\V_{\alpha,\beta}\|_2^2\le\frac C{n^{\beta-\alpha+1}}\sum_{j=1}^n\sum_{j_1=1}^{n}
\sum_{j_2=1}^{n}\dots \sum_{j_{\beta-\alpha}=1}^{n}\sum_{k=1}^{n}
\E [X^{(\alpha)}_{j,j_1}X^{(\alpha+1)}_{j_1,j_2}\cdots X^{(\beta)}_{j_{\beta-\alpha},k}]^2
\end{equation*}
By independents of random variables, we get
\begin{equation*}
\E\|\V_{\alpha,\beta}\|_2^2\le  Cn
\end{equation*}
Thus the Lemma is proved.
\end{proof}

\begin{lemma}\label{norm4}
Under the condition of Theorem \ref{th:main} and assumption (\ref{as1}) we have, for any $j,k=1,\ldots,n $, and $r\ge1$,
\begin{equation}\label{new1}
\E\|\V_{\alpha,\beta}\mathbf e_k\|_2^{2r}\le C_r,\quad \E\|\V_{\alpha,\beta}\mathbf e_{j+n}\|_2^{2r}\le C_r
\end{equation}
and
\begin{equation}\label{new2}
\E\|{\ee_j}^T\V_{\alpha,\beta}\|_2^{2r}\le C_r,\quad \E\|{\ee_{k+n}}^T\V_{\alpha,\beta}\|_2^{2r}\le C_r,
\end{equation}
with some positive constant $C_r$ depending on $r$. Moreover, for any $q=1,\ldots,m$ and any $l,s=1,\ldots,n$,
\begin{equation}\label{new3}
\E\Big\{\|{\ee_j}^T\V_{\alpha,\beta}\|_2^{2r}\Big|X^{(q)}_{ls},X^{(q)}_{sl}\Big\}\le C_r.
\end{equation}
and
\begin{equation}\label{new4}
\E\Big\{\|\V_{\alpha,\beta}\mathbf e_{j+n}^{(\beta)}\|_2^{2r}\Big|X^{(\nu)}_{lq},X^{(q)}_{sl}\Big\}\le C_r.
\end{equation}
\end{lemma}
\begin{proof}
By definition  of the  matrices $\V_{\alpha,\beta}$, we may write
\begin{equation*}
\|\ee_j^T\V_{\alpha,\beta}\|_2^{2}=\frac1{n^{\beta-\alpha}}\sum_{l=1}^{n}\left|\sum_{j_{\alpha}=1}^{n}\cdots\sum_{j_{\beta-1}=1}^{n}
X_{jj_{\alpha}}^{(\alpha)}\cdots X_{j_{\beta-1}l}^{(\beta)}\right|^2
\end{equation*}
Using this representation, we get
\begin{align}\label{mom10}
&\E\|\ee_j^T\V_{\alpha,\beta}\|_2^{2r}=\frac1{n^{(\beta-\alpha)r}}\nonumber\\
&\times\sum_{l_1=1}^{n}\cdots\sum_{l_r=1}^{n}\E\prod_{q=1}^r
\left(\sum_{j_{\alpha}=1}^{n}\cdots\sum_{j_{\beta-1}=1}^{n}
\sum_{\widehat j_{\alpha}=1}^{n}\cdots
\sum_{\widehat j_{\beta-1}=1}^{n}A^{(l_q)}_{(j_{\alpha},\ldots,j_{\beta-1},\widehat j_{\alpha},\ldots,
\widehat j_{\beta-1})}\right)
\end{align}
where
\begin{align}\label{per}
&A^{(l_q)}_{(j_{\alpha},\ldots,j_{\beta-1},\widehat j_1,\ldots,\widehat j_{\beta-1})}:=\nonumber\\
&\qquad\qquad X_{jj_{\alpha}}^{(\alpha)}
\overline X_{j \widehat j_{\alpha}}^{(\alpha)}X_{j_{\alpha}j_{\alpha+1}}^{(\alpha+1)}
 X_{{\widehat j}_{\alpha} \widehat j_{\alpha+1}}^{(\alpha+1)}\cdots X_{j_{\beta-2}j_{\beta-1}}^{(\beta-1)}
 X_{\widehat j_{\beta-2}\widehat  j_{\beta-1}}^{(\beta-1)}X_{j_{\beta-1}l_q}^{(\beta)}
 X_{\widehat j_{\beta-1}l_q}^{(\beta)}.
\end{align}
Rewriting the product on  the r.h.s of (\ref{mom10}), we get
\begin{align}\label{per1}
 \E\|\mathbf e_j^T\mathbf V_{\alpha,\beta}\|_2^{2r}=\frac1{n^{(\beta-\alpha)r}}
{\sum}^{**}\E\prod_{q=1}^rA^{(l_q)}_{(j_{\alpha}^{(q)},\ldots,j_{\beta-1}^{(q)},
{\widehat j}_1^{(\nu)},\ldots,{\widehat j}_{\beta-1}^{(q)})},
\end{align}
where ${\sum}^{**}$ is taken  over all set of indices $j_{\alpha}^{(q)},\ldots, j_{\beta-1}^{(q)}, l_q$ and
${\widehat j}_{\alpha}^{(\nu)},\ldots,{\widehat j}_{\beta-1}^{(q)}$ where
$j_k^{(q)},{\widehat j}_k^{(q)}=1,\ldots,p_k$, $k=\alpha,\ldots,\beta-1$, $l_q=1,\ldots,n$ and  $q=1,\ldots,r$.
Note that the summands in the right hand side of (\ref{per}) is equal 0 if there is at least one term in the product (\ref{per})
which appears only one time. This implies that  the summands in the  right hand side of (\ref{per1}) is not equal zero only if the union of all sets of indices in r.h.s of (\ref{per}) consist from at least $r$ different terms and each term appears at least twice.

Introduce the random variables, for $q=\alpha+1,\ldots, \beta-1$,
\begin{align*}
&\zeta^{(q)}_{j^{(1)}_{q-1},\ldots,j^{(r)}_{q-1},j^{(1)}_{q},\ldots,j^{(r)}_{q},
{\widehat j}^{(1)}_{q-1},\ldots,{\widehat j}^{(r)}_{q-1},{\widehat j}^{(1)}_{q},\ldots,
{\widehat j}^{(r)}_{q}}:= \\
&\qquad\qquad\qquad\qquad\qquad\qquad X^{(q)}_{j^{(1)}_{q-1},j^{(1)}_{q}}
\cdots X^{(q)}_{j^{(r)}_{q-1},j^{(r)}_{q}}
{ X}^{(q)}_{{\widehat j}^{(1)}_{q-1},{\widehat j}^{(1)}_{q}},
\cdots { X}^{(q)}_{{\widehat j}^{(r)}_{q-1},{\widehat j}^{(r)}_{q}},
\end{align*}
and
\begin{align*}
\zeta^{(\alpha)}_{j^{(1)}_{1},\ldots,j^{(r)}_{1},
{\widehat j}^{(1)}_{1},\ldots,{\widehat j}^{(r)}_{1}}&:=X^{(\alpha)}_{jj_1^{(\alpha)}}
\cdots X^{(\alpha)}_{j^{(r)}_{a}j^{(r)}_{a+1}}
{ X}^{(\alpha)}_{j{\widehat j}^{(1)}_{a}}
\cdots { X}^{(\alpha)}_{{\widehat j}^{(r)}_{a},{\widehat j}^{(r)}_{a+1}}\\
\zeta^{(\beta)}_{j^{(1)}_{\beta-1},\ldots,j^{(r)}_{\beta-1},
{\widehat j}^{(1)}_{\beta-1},\ldots,{\widehat j}^{(r)}_{\beta-1},l_q}
&:=X^{(\beta)}_{j^{(1)}_{\beta-1}j^{(1)}_{\beta}}
\cdots X^{(\beta)}_{j^{(r)}_{\beta-1}l_q}
{ X}^{(\beta)}_{{\widehat j}^{(1)}_{\beta-1},l_q},
\cdots { X}^{(\beta)}_{{\widehat j}^{(r)}_{\beta-1},l_q}.
\end{align*}
 Assume that the set of indices $j^{(1)}_{\alpha},\ldots,j^{(r)}_{\alpha},
{\widehat j}^{(1)}_{\alpha},\ldots,{\widehat j}^{(r)}_{\alpha}$ contains $t_{\alpha}$ different indexes, say
$i_1^{(\alpha)},\ldots,i_{t_{\alpha}}^{(\alpha)}$ with multiplicities $k_1^{(\alpha)},\ldots,k_{t_{\alpha}}^{(\alpha)}$
respectively, $k_1^{(\alpha)}+\ldots+k_{t_{\alpha}}^{(\alpha)}=2r$.
Note that $\min\{k_1^{(\alpha)},\ldots,k_{t_{\alpha}}^{(\alpha)}\}\ge 2$.
Otherwise,\newline $|\E\zeta^{(\alpha)}_{j^{(1)}_{a},\ldots,j^{(r)}_{a},
{\widehat j}^{(1)}_{\alpha},\ldots,{\widehat j}^{(r)}_{\alpha}}|=0$. By assumption (\ref{as1}), we have
\begin{equation}\label{mom1}
|\E\zeta^{(\alpha)}_{j^{(1)}_{\alpha},\ldots,j^{(r)}_{\alpha},
{\widehat j}^{(1)}_{\alpha},\ldots,{\widehat j}^{(r)}_{\alpha}}|\le C(\tau_n\sqrt n)^{2r-2t_{\alpha}}
\end{equation}
Similar bounds we get for $|\E
\zeta^{(\beta)}_{j^{(1)}_{\beta-1},\ldots,j^{(r)}_{1},
{\widehat j}^{(1)}_{\beta-1},\ldots,{\widehat
  j}^{(r)}_{\beta-1},l_q}|$. Assume that
the set of indexes $\{j^{(1)}_{\beta-1},\ldots,j^{(r)}_{\beta-1}$,
${\widehat j}^{(1)}_{\beta-1},\ldots,{\widehat j}^{(r)}_{\beta-1}\}$ contains $t_{\beta-1}$ different indices, say,
$i_1^{(\beta-1)},\ldots,i_{t_{\beta-1}}^{(\alpha)}$ with multiplicities\newline
$k_1^{(\beta-1)},\ldots,k_{t_{\beta-1}}^{(\alpha)}$
respectively, $k_1^{(\beta-1)}+\ldots+k_{t_{\beta-1}}^{(\alpha)}=2r$. Then
\begin{equation}\label{mom2}
|\E
\zeta^{(\beta)}_{j^{(1)}_{\beta-1},\ldots,j^{(r)}_{1},
{\widehat j}^{(1)}_{\beta-1},\ldots,{\widehat j}^{(r)}_{\beta-1},l_q}|\le C(\tau_n\sqrt n)^{2r-2t_{\beta-1}}
\end{equation}
Furthermore, assume that for $\alpha+1\le q\le \beta-2$ there are $t_{q}$
{different  pairs of indices}, say, $(i_{\alpha},
i'_{\alpha}),
\ldots(i_{t_{\beta}},i'_{t_{\beta}})$ in the set\newline $\{j^{(1)}_{\alpha},\ldots,j^{(r)}_{\alpha},
{\widehat j}^{(1)}_{\alpha},\ldots,{\widehat j}^{(r)}_{\alpha},\ldots,j^{(1)}_{\beta-1},\ldots,j^{(r)}_{\beta-1},
{\widehat j}^{(1)}_{\beta-1},\ldots,{\widehat j}^{(r)}_{\beta-1},l_1,l_r\}$
with multiplicities\newline$k_1^{(q)},\ldots,k_{t_{q}}^{(q)}$.
Note that
\begin{equation*}
 k_1^{(q)}+\ldots+k_{t_{q}}^{(q)}=2r
\end{equation*}
and
\begin{equation}\label{mom3}
|\E\zeta^{(q)}_{j^{(1)}_{q-1},\ldots,j^{(r)}_{q-1},j^{(1)}_{q},\ldots,j^{(r)}_{q},
{\widehat j}^{(1)}_{q-1},\ldots,{\widehat j}^{(r)}_{q-1},{\widehat j}^{(1)}_{q},\ldots,{\widehat j}^{(r)}_{q}}|\le
C(\tau_n\sqrt n)^{2r-2t_{q}}.
\end{equation}
Inequalities (\ref{mom1})-(\ref{mom3}) together yield
\begin{equation}\label{mom11}
 |\E\prod_{q=1}^rA^{(l_q)}_{(j_{\alpha}^{(q)},\ldots,
j_{\beta-1}^{(q)},{\widehat j}_1^{(q)},\ldots,{\widehat j}_{\beta-1}^{(q)})}|
\le C(\tau_n\sqrt n)^{2r(\beta-\alpha)-2(t_1+\ldots+t_{\beta-\alpha})}.
\end{equation}
It is straightforward to check that the  number $\mathcal N(t_{\alpha},\ldots,t_{\beta})$ of sequences of indices
\newline$\{j^{(1)}_{\alpha},\ldots,j^{(r)}_{\alpha},
{\widehat j}^{(1)}_{\alpha},\ldots,{\widehat j}^{(r)}_{\alpha},\ldots,j^{(1)}_{\beta-1},\ldots,j^{(r)}_{\beta-1},
{\widehat j}^{(1)}_{\beta-1},\ldots,{\widehat j}^{(r)}_{\beta-1},l_1,\ldots,l_r\}$ with $t_{\alpha},\ldots,t_{\beta}$
{of different pairs}  satisfies the inequality
\begin{equation}\label{finish}
 \mathcal N(t_{\alpha},\ldots,t_{\beta})\le Cn^{t_{\alpha}+\ldots+t_{\beta}},
\end{equation}
with $1\le t_i\le r,\quad i=\alpha,\ldots,\beta$. Note that in the case $t_{\alpha}=\cdots=t_b=r$ the inequalities~\eqref{mom1}--~\eqref{mom3} imply
\begin{equation}\label{mom4}
 \E\zeta^{(q)}_{j^{(1)}_{q-1},\ldots,j^{(r)}_{q-1},j^{(1)}_{\nu},\ldots,j^{(r)}_{q},
{\widehat j}^{(1)}_{q-1},\ldots,{\widehat j}^{(r)}_{q-1},{\widehat j}^{(1)}_{\nu},\ldots,{\widehat j}^{(r)}_{q}}\le C
\end{equation}
The inequalities~\eqref{finish},~\eqref{mom11},~\eqref{mom4}, and the representation~\eqref{mom10} together conclude the proof of inequalities~\eqref{new1} and~\eqref{new2}. To prove the inequalities~\eqref{new3}, ~\eqref{new4} note that in the case $q\notin[\alpha,\beta]$ and $m-q\notin[\alpha,\beta]$ we have
\begin{align*}
\E\Big\{\|{\ee_j}^T\V_{\alpha,\beta}\|_2^{2r}\Big|X^{(q)}_{ls},X^{(q)}_{ls}\Big\}=\E\|{\ee_j}^T\V_{\alpha,\beta}\|_2^{2r}\\
\E\Big\{\|\V_{\alpha,\beta}\ee_{j+n}\|_2^{2r}\Big|X^{(q)}_{ls},X^{(q)}_{sl}\Big\}=\E\|\V_{\alpha,\beta}\ee_{j+n}\|_2^{2r}.
\end{align*}
Thus in the case $q\notin[\alpha,\beta]$ and $m-q\notin[\alpha,\beta]$ the inequalities~\eqref{new3} and~\eqref{new4} are proved. Consider now the case $q\in[\alpha,\beta]$ and $m-q\notin[\alpha,\beta]$. In this case we may write
\begin{align}\label{n1}
\V_{\alpha,\beta}=\V_{\alpha,q-1}(\mathbf H^{(q,l,s)}+X^{(q)}_{ls}{\mathbf e_l}{\ee_q}^T+X^{(q)}_{sl}{\ee_s}
{\ee_l}^T)\V_{q+1,\beta},
\end{align}
where the matrix $\mathbf H^{(q,l,s)}$ is  obtained from the matrix $\HH^{(q)}$ by replacement the entries $X_{ls}^{(q)}$ and $X_{sl}^{(q)}$ by zero. Note that the matrix $\HH^{(q,l,s)}$ and random variables $X_{ls}^{(q)}$ and $X_{ls}^{(q)}$ are independent. Let $\V_{\alpha,\beta}^{(q,l,s)}=\V_{\alpha,q-1}\HH^{(q,l,s)}\mathbf V_{q+1,\beta}$. We may rewrite~\eqref{n1} in the form
\begin{equation}\label{n2}
\V_{\alpha,\beta}=\V_{\alpha,\beta}^{(\nu,l,q)}+\frac1{\sqrt n}X^{(q)}_{ls}\V_{\alpha,\nu-1}{\ee_l}
{\ee_s}^T\V_{q+1,\beta}+\frac1{\sqrt n}X^{(q)}_{sl}\V_{\alpha,q-1}{\ee_s}
{\ee_l}^T\V_{q+1,\beta}
\end{equation}
From the independence of $\V_{\alpha,q-1}$, $\V_{q+1,\beta}$, $X^{(q)}_{ls}$, $X^{(q)}_{sl}$ and $|X^{(q)}_{ls}|/\sqrt n\le \tau_n$, the equality~\eqref{n2} it follows that
\begin{align*}
\E\Big\{\|\V_{\alpha,\beta}\ee^{(q)}_j\|_2^{2r}\Big|\xi^{(q)}_
{ls},\xi^{(q)}_{ls}\Big\}\le& 2^r\Big(\E\|\V_{\alpha,\beta}^{(q,l,s)}\ee_j\|_2^{2r}\\
&+\tau_n\E\|\V_{\alpha,\nu-1}{\ee_l}\|_2^{2r}\E\|{\ee_q}^T\V_{q+1,\beta}\ee_j\|_2^{2r}\Big).
\end{align*}
The last inequality concludes the proof of inequality~\eqref{new3} in the case $q\in[\alpha,\beta]$ and $m-q\notin[\alpha,\beta]$. The proof of inequality~\eqref{new4} is similar. The proof of both inequalities (\ref{new3}) and~\eqref{new4} in the cases $q\notin[\alpha,\beta]$ and $m-q\in[\alpha,\beta]$ and $q\in[\alpha,\beta]$ and $m-q\in[\alpha,\beta]$ is analogously. Thus Lemma~\ref{norm4} is proved.
\end{proof}

\begin{lemma}\label{var0}
Under conditions of Theorem \ref{th:main} assuming (\ref{as1}), we have
\begin{equation*}
\E\left|\frac1n \left(\Tr \RR-\E\Tr \RR \right )\right|^2\le \frac C{nv^2}.
\end{equation*}
\end{lemma}
\begin{proof}
We define the following matrices
\begin{equation*}
\HH^{(q,j)}=\HH^{(q)}-\ee_j\ee_j^T\HH^{(q)}-\HH^{(q)}\ee_j\ee_j^T,
\end{equation*}
and
\begin{equation*}
{\widetilde{\HH}}^{(m-q+1,j)}={{\HH}}^{(m-q+1)}-{{\HH}}^{(m-q+1)}\ee_{j+n}\ee_{j+n}^T-\ee_{j+n}\ee_{j+n}^T{{\HH}}^{(m-q+1)},
\end{equation*}
for $q=1,\ldots,m$ and $j=1,\ldots,n$. For simplicity we shall assume that  $q\le m-q+1$.
Define 
$$
\V^{(q,j)}=\prod_{\beta=1}^{q-1}\HH^{(\beta)}\,\HH^{(q,j)}\prod_{\beta=q+1}^{m-q}\mathbf H^{(\beta)}{\widetilde{\HH}}^{(m-q+1,j)}\prod_{\beta=m-q+2}^{m}\HH^{(\beta)}.
$$
Let $\V^{(q,j)}(z)=\V^{(q,j)}\J-\J(z)$.  We shall use the following inequality. For any Hermitian matrices  $\mathbf A$ and $\mathbf B$ with spectral distribution function $F_A(x)$ and $F_B(x)$ respectively, we have
\begin{equation}\label{trace}
|\Tr (\mathbf A-\alpha\mathbf I)^{-1}-\Tr (\mathbf B-\alpha\mathbf I)^{-1}|\le
\frac {\text{\rm rank}(\mathbf A-\mathbf B)}{v},
\end{equation}
where $\alpha=u+iv$. It is straightforward to show that
\begin{equation}\label{rank}
\text{\rm rank}(\mathbf V(z)-\mathbf V^{(q,j)}(z))=\text{\rm rank}(\mathbf V\mathbf J-\mathbf V^{(q,j)}\mathbf J)\le 4m.
\end{equation}
The inequalities~\eqref{trace} and~\eqref{rank} together imply
\begin{equation*}
\left |\frac1{2n}(\Tr \RR-\Tr \RR^{(q,j)})\right|\le \frac C{nv}.
\end{equation*}

After this remark we may apply a standard martingale expansion procedure. We  introduce
$\sigma$-algebras $\mathcal F_{q,j}=\sigma\{X^{(q)}_{lk},\,j< l,k\le n; X^{(\beta)}_{ps}$,
$\beta=q+1,\ldots m, \,p,s=1,\ldots,n,\}$ and use the representation
\begin{equation*}
\Tr\RR-\E\Tr\RR=\sum_{q=1}^m\sum_{j=1}^{n}(\E_{q,j-1}\Tr\RR-\E_{q,j}\Tr\RR),
\end{equation*}
where $\E_{q,j}$ denotes  conditional expectation given the  $\sigma$-algebra $\mathcal F_{q,j}$. Note that $\mathcal F_{q,n}=\mathcal F_{q+1,0}$ and $\E_{q,j-1}\Tr\mathbf R^{(q,j)}=\E_{q,j}\Tr\mathbf R^{(q,j)}$.
\end{proof}
%%%%%%%%%%%%%%%%%%%%%%%%%%%%%%%%%%%%%%%%%%%%%%%%%%%%%%%%%%%%%%%%%%%%%%%%%%%%%%%%%%%%%%%%%%%%%%%%%%%%%%%%%%%%%%%%%%%%%%%%%%%%%%%%%%%%%%%%%%%%%%%%%%
\begin{lemma}\label{var1}
Under the conditions of Theorem \ref{th:main} we have, for $1\le a\le m$,
\begin{equation*}
\E\left|\frac1n\left(\sum_{k=1}^{n}[\V_{a+1,m}\J\RR\V_{1,m-a}]_{k,k+n}-\E\sum_{j=1}^{n}[\V_{a+1,m}\J\RR\V_{1,m-a}]_{kk+{n}}\right)\right|^2\le \frac C{n v^4}.
\end{equation*}
and, for $1\le a\le m-1$,
\begin{equation*}
\E\left|\frac1n\left(\sum_{k=1}^{n}[\V_{m-a+2,m}\J\RR\V_{1,m-a+1}]_{k,k}-\E\sum_{j=1}^{n}[\V_{m-a+2,m}\J\RR\V_{1,m-a+1}]_{kk}\right)\right|^2
\le \frac C{n v^4}.
\end{equation*}
\end{lemma}
%%%%%%%%%%%%%%%%%%%%%%%%%%%%%%%%%%%%%%%%%%%%%%%%%%%%%%%%%%%%%%%%%%%%%%%%%%%%%%%%%%%%%%%%%%%%%%%%%%%%%%%%%%%%%%%%%%%%%%%%%%%%%%%%%%%%%%%%%%%%%%%%%
\begin{proof}
We prove the first inequality only. The proof of the other one is  similar. Let $\HH^{(q,j)}$ and
${\widetilde{\HH}}^{(m-q+1,j)}$ be the matrices defined in the previous Lemma, for $q=1,\ldots,m$ and for $j=1,\ldots,n$. We introduce as well the matrices $\mathbf X^{(q,j)}=\mathbf X^{(q)}-\ee_j\ee_j^T\mathbf X^{(q)}-\mathbf X^{(q)}\ee_j\ee_j^T$. Note that the matrix $\mathbf X^{(q,j)}$ is obtained from the matrix $\mathbf X^{(q)}$ by replacing its $j$-th row and $j$th column by a row  and column  of zeros. Similar to the proof  of the
previous Lemma we introduce the matrices $\mathbf V^{(q,j)}_{c,d}$ by replacing in the definition of
$\V_{c,d}$ the  matrix $\HH^{(q)}$ by $\HH^{(q,j)}$ and the matrix $\HH^{(m-q+1)}$ by ${\widetilde{\HH}}^{(m-q+1,j)}$. For instance, if $c\leq m-q+1\le d$ we get
\begin{equation*}
\V^{(q,j)}_{c,d}=\prod_{\beta=c}^{q-1}\HH^{(\beta)}\,\HH^{(q,j)}\prod_{\beta=q+1}^{m-q}\HH^{(\beta)}{\widetilde{\HH}}^{(m-q+1,j)}\prod_{\beta=m-q+1}^{d}\HH^{(\beta)}.
\end{equation*}
Let $\V^{(q,j)}:= \V_{1,m}^{(q,j)}$ and $\RR^{(j)}:= (\V^{(q,j)}(z)-\alpha\mathbf I)^{-1}$. Introduce  the
following quantities, for $q=1\ldots,m$ and $j=1,\ldots,n$,
\begin{equation*}
\Xi_{q,j}:=\sum_{k=1}^n[\V_{a+1,m}\J\RR\V_{1,m-a+1}]_{kk+n}- \sum_{k=1}^n[\V^{(q ,j)}_{a+1,m} \J\RR^{(q,j)}\V^{(q,j)}_{1,m-a+1}]_{kk+n}
\end{equation*}
We represent them  in the following form
\begin{equation*}
\Xi_{q,j}:= \Xi_{q,j}^{(1)}+\Xi_{q,j}^{(2)}+\Xi_{q,j}^{(3)},
\end{equation*}
where
\begin{align*}
\Xi_{q,j}^{(1)}&=\sum_{k=1}^{n}[(\V_{a+1,m}-\V_{a+1,m}^{(q,j)})\J\RR\V_{1,m-a+1}]_{k,k+n},\\
\Xi_{q,j}^{(2)}&= \sum_{k=1}^{n}[\V_{a+1,m}^{(q,j)}\J(\RR-\RR^{(q,j)})\J\V_{1,m-a+1}]_{k,k+n},\\
\Xi_{q,j}^{(3)}&= \sum_{k=1}^{n}[\V^{(j)}_{a+1,m}\J\RR^{(q,j)}(\ V_{1,m-a+1}-\V_{1,m-a+1}^{(q,j)})]_{k,k+n}.
\end{align*}
Note that
\begin{align*}
\V_{a+1,m}-\mathbf V^{(q,j)}_{a+1,m}&=\V_{a+1,q-1}(\HH^{(q)}-\HH^{(q,j)})\mathbf V_{q+1,m}\\
&+\V_{a+1,q-1}\HH^{(q,j)}\V_{q+1,m-\nu}(\widetilde {\HH}_{m-q+1}-{\widetilde {\HH}}_{m-q+1}^{q,j})\V_{m-q+2,m}.
\end{align*}
By definition of the matrices $\HH^{q,j}$ and ${\widetilde{\HH}}^{m-q+1,j}$, we have
\begin{align*}
\sum_{k=1}^{n}[(\V_{a+1,m} -\V_{a+1,m}^{(q,j)})\J\RR\V_{1,m-q+1}]_{k,k+n}=
[\V_{q+1,m}\J\RR\V_{1,m-a+1}\mathbf{\widetilde J} \V_{a+1,q}]_{j,j}\\
+[\V_{m-q+2,m}\J\RR\V_{1,m-a+1}\mathbf{\widetilde J}\V_{a+1,m-a+1}]_{j+n,j+n},
\end{align*}
where
$$
\mathbf{\widetilde J}=\left(\begin{matrix}{\mathbf O\quad\mathbf
I}\\{\mathbf O\quad\mathbf O}\end{matrix}\right)
$$

This equality implies that
\begin{align*}
|\Xi_{q,j}^{(1)}|&\le |[\V_{q+1,m}\J\RR\V_{1,m-a+1}\mathbf{\widetilde J} \V_{a+1,q}]_{j,j+n}|\\
&\qquad\qquad+|[\V_{m-q+2,m}\J\RR\V_{1,m-a+1}\mathbf{\widetilde J}\V_{a+1,m-q+1}]_{j+n,j+n}|.
\end{align*}
Using the obvious  inequality $\sum_{j=1}^n a_{jj}^2\le \|\mathbf A\|_2^2$ for any matrix
$\mathbf A=(a_{jk})$, $j,k=1,\ldots,n$, we get
\begin{align*}
T_1:=\sum_{j=1}^n\E|\Xi_j^{(1)}|^2\le &\E\|\V_{q+1,m}\J\RR\V_{1,m-a+1} \mathbf{\widetilde J}\V_{a+1,q}\|_2^2\notag\\&+\E\|\V_{m-q+2,m}\J\RR\V_{1,m-a+1}\mathbf{\widetilde J}\V_{a+1,m-q+1}\|_2^2.
\end{align*}
By Lemma \ref{norm2}, we get
\begin{equation}\label{T1}
T_1\le \frac{C}{v^2}\E\|\mathbf V_{a+1,m}\mathbf V_{1,m-a+1}\|_2^2\le \frac {Cn}{v^2}
\end{equation}
Consider now the term
\begin{equation*}
 T_2=\sum_{j=1}^n\E|\Xi_{q,j}^{(2)}|^2.
\end{equation*}
Using that $\mathbf R-\mathbf R^{(j)}=-\mathbf R^{(j)}(\mathbf V(z)-\mathbf
V^{(q,j)}(z))\mathbf R$, we get
\begin{align*}
|\Xi_{q,j}^{(2)}|&\le |\sum_{k=1}^{n}[\V^{(q,j)}_{a,m}\J\RR\V_{1,q-1}\ee_j\ee_j^T \V_{q,m}
\mathbf R\V_{1,b}]_{k,k+n}|\\
&\qquad\le [\J\HH^{(\alpha+1)}\V_{\alpha+2,m-\alpha}\HH^{(m-\alpha+1,j)}\V_{m-\alpha+2,m}\RR\V_{1,m-\alpha}\V^{(j)}_{\alpha+1,m}
\J\RR\V_{1,\alpha}]_{jj}.
\end{align*}
This implies that
\begin{equation*}
 T_{2}\le C\E\|[\V_{q+1,m}\J\RR\V_{1,b}\V_{a,m}\J\RR\V_{1,q}\|_2^2.
\end{equation*}
It is straightforward to check that
\begin{equation}\label{t2}
 T_{2}\le \frac C{v^4}\E\|\mathbf V_{1,\alpha}\mathbf J\mathbf H^{(\alpha+1)}
\mathbf V_{\alpha+2,m-\alpha}\mathbf H^{(m-\alpha+1,j)}
\mathbf V_{m-\alpha+2,m}\|_2^2=\E\|\mathbf Q\|_2^2
\end{equation}
The matrix on the right hand side of equation (\ref{t2}) may be  represented in the following form
\begin{equation*}
 Q=\prod_{q=1}^m{\HH^{(q)}}^{\varkappa_{q}},
\end{equation*}
where $\varkappa_{q}=0$ or $\varkappa_{q}=1$ or $\varkappa_{q}=2$. Since $X^{(q)}_{ss}=0$, for $\varkappa=1$ or $\varkappa=2$, we have
\begin{equation*}
\E|{\HH^{(q)}}^{\varkappa}_{kl}|^2\le \frac C{n}.
\end{equation*}
This implies that
\begin{equation}\label{T2}
T_2\le Cn.
\end{equation}
Similar we prove that
\begin{equation}\label{T3}
 T_3:=\sum_{j=1}^n\E|\Xi_{q,j}^{(3)}|^2\le Cn.
\end{equation}
The inequalities~\eqref{T1},~\eqref{T2} and~\eqref{T3} together imply
\begin{equation*}
 \sum_{j=1}^n\E|\Xi_{q,j}|^2\le Cn
\end{equation*}
Applying now a martingale expansion with respect to the $\sigma$-algebras $\mathcal F_j$ generated by the random variables $X_{kl}^{(\alpha+1)}$ with $1\le k\le j$, $1\le l\le n$ and all other random variables $X^{(q)}_{sl}$
except $q=\alpha+1$, we get
\begin{equation*}
\E\left|\frac1n\left(\sum_{k =1}^n[\V_{\alpha+1,m}\J\RR
\V_{1,m-\alpha}]_{kk+n}-\E\sum_{j=1}^n[\V_{\alpha+1,m}\J\RR \V_{1,m-\alpha}]_{kk+n}\right)\right|^2\le \frac C{n v^4}.
\end{equation*}
Thus the Lemma is proved.
\end{proof}

%\newpage

\bibliographystyle{plain}
\bibliography{literatur}
\end{document}